\documentclass[a4paper,11pt]{amsart}
\usepackage{amsmath,amssymb,latexsym}
\usepackage{enumitem}
\usepackage{color}
\usepackage{dsfont}
\usepackage{stmaryrd}
\catcode`\é = \active\def é{\'e} \catcode`\è = \active\def è{\`e}
\catcode`\ç = \active\def ç{\c{c}} \catcode`\à = \active\def
à{\`a} \catcode`\ê = \active \def ê{\^e} \catcode`\ù = \active
\def ù{\`u} \catcode`\ô = \active \def ô{\^o} \catcode`\û =
\active \def û{\^u} \catcode`\î = \active \def î{\^{\i}}
\catcode`\â = \active \def â{\^a} \catcode`\ö = \active \def
ö{{\"o}}

\newtheorem{thm}{Theorem}[section]  \newtheorem{cor}[thm]{Corollary}
\newtheorem{lem}[thm]{Lemma}   \newtheorem{defn}[thm]{Definition}
  \newtheorem{prop}[thm]{Proposition}

\def\remark{\refstepcounter{thm}\bigskip\noindent\bf Remark \thethm\rm\ }
\newcommand{\preuve}[1][\!\!]{\bigskip\noindent{\bf Proof #1. \ \ }}
\def\fin{\hfill$\Box$\\}

 \def\hhh{{\mathcal H}}

\def\sss{{\mathcal S}}

\def\Bk{\color{black}}

\def\R{\mathbb R}\def\N{\mathbb N}\def\Z{\mathbb Z}
  \def\S{\mathbb S} \def\T{\mathbb T}   
\def\D{\partial}\def\eps{\varepsilon}\def\phi{\varphi}
\def\d{{\rm d}}
\def\norm#1{\left\Vert#1\right\Vert}

\def\abs#1{\left\vert#1\right\vert}
\def\set#1{\left\{#1\right\}}
\def\seq#1{\left<#1\right>}
\def\sep#1{\left(#1\right)}\def\adf#1{\left[#1\right]}
\def\Re{{\mathrm Re\,}}

\def\wick{{\rm Wick}}
\def\seta{\seq{\eta}}
\def\sxi{\seq{\xi}}

\newcommand{\dv}{{{\rm d}v}}

\def\exp{\textrm{e}}
\newcommand{\ddt}{{\frac{{{\rm d}}}{{{\rm d}t}}}}
\def\hf{\widehat{f}}

\date{January 30, 2006}
\title[Short time diffusion with fractional collision kernel]{Short time diffusion properties of inhomogeneous kinetic equations with fractional collision kernel}

\author{{\sc Fr\'{e}d\'{e}ric H\'{e}rau}}
\address{Laboratoire de Math\'{e}matiques Jean Leray, CNRS,
Universit\'{e} de Nantes, 44322 Nantes, France. E-mail: {\tt herau@univ-nantes.fr}}

\author{{\sc Daniela Tonon}}
\address{CEREMADE, UMR 7534, Universit\'e Paris IX-Dauphine, PSL Research University,
Place du Mar\'echal de Lattre de Tassigny, 75775 Paris Cedex 16, France.
E-mail: {\tt tonon@ceremade.dauphine.fr }}

\author{{\sc Isabelle Tristani}}
\address{D\'epartement de mathématiques et applications, \'Ecole normale supérieure, CNRS, PSL Research University, 45 rue d'Ulm, 75005 Paris, France
E-mail: {\tt isabelle.tristani@ens.fr}}

\begin{document}

\date{\today}

 \maketitle

 \thispagestyle{empty}
 \begin{abstract} We prove regularization properties in short time for inhomogeneous kinetic equations whose collision kernel behaves like a fractional power of the Laplacian in velocity. We treat a fractional Kolmogorov equation and the linearized Boltzmann equation without cutoff (for hard potentials).
\end{abstract}

\tableofcontents

\section{Introduction and results}

\subsection{Models}
In this paper, we consider two kinetic inhomogeneous equations
on $\R^+\times\T^d_x \times \R^d_v$ ($\T^d$ being the $d$-dimensional torus) with collision kernel having the behavior of a fractional Laplacian in velocity:

\begin{equation}
  \label{FKE}
  \left\{
  \begin{aligned}
   &  \D_{t}f+v\cdot\nabla_{x}f = { Lf } \\
    &   f|_{t=0} = f_0,
    \end{aligned}
    \right.
\end{equation}
where $f = f(t,x,v)$ is the distribution of particles, and $L$ is the collision kernel, roughly behaving like a fractional power of the Laplacian in velocity, and acting only in velocity.

\subsubsection{Fractional Kolmogorov equation} The simplest model entering in this family is the {fractional Kolmogorov equation}: for $s\in (0,1]$, the corresponding collision kernel is given by
\begin{equation} \label{FKO}
L := -( 1 -\Delta_v)^{s}.
\end{equation}

\subsubsection{Linearized Boltzmann equation} We also deal with a more complicated model associated to the {linearized Boltzmann operator without cutoff for hard potentials} in dimension $d=3$.
Let us describe it more precisely.
The Boltzmann collision operator is defined as
	\begin{equation} \label{eq:Q}
	Q(g,f):=\int_{\R^3 \times \mathbb{S}^2} B(v-v_*,\sigma) \left[g'_*f' - g_* f \right] \, \d \sigma \, \d v_*.
	\end{equation}
Here and below, we are using the shorthand notations $f=f(v)$, $g_*=g(v_*)$, $f'=f(v')$ and $g'_*=g(v'_*)$.
In this expression, $v$, $v_*$ and $v'$, $v'_*$ are the velocities of a pair of particles before and after collision.
We make a choice of parametrization of the set of solutions to the conservation of momentum and energy (physical laws of elastic collisions):
	$$
	v+v_*=v'+v'_* ,
	$$
	$$
	|v|^2+|v_*|^2=|v'|^2+ |v'_*|^2,
	$$
so that the post-collisional velocities are given by:
	$$
	v'=\frac{v+v_*}{2} + \frac{|v-v_*|}{2} \sigma, \quad v'_*=\frac{v+v_*}{2} - \frac{|v-v_*|}{2} \sigma, \quad \sigma \in \mathbb{S}^2.
	$$
The Boltzmann collision kernel $B(v-v_*,\sigma)$ only depends on the relative velocity $|v-v_*|$ and on the deviation angle $\theta$ through $\cos \theta = \langle \kappa, \sigma \rangle$
where $\kappa = (v-v_*)/|v-v_*|$ and $\langle \cdot, \cdot \rangle$ is the usual scalar product in $\R^3$.
By a symmetry argument, one can always reduce to the case where $B(v-v_*, \sigma)$ is supported on $ \langle \kappa, \sigma \rangle \geq 0$ i.e. $0 \leq \theta \leq \pi/2$.
So, without loss of generality, we make this assumption.
In the sequel, we shall be concerned with the case when the kernel $B$ satisfies the following conditions (which include the physical case of the so-called hard potentials):
	\begin{itemize}
	\item it takes product form in its arguments as
		\begin{equation} \label{eq:B}
			B(v-v_*,\sigma) = \Phi(|v-v_*|) \, b(\cos \theta);
		\end{equation}
	\item the angular function $b$ is locally smooth, and has a nonintegrable singularity for $\theta \rightarrow 0$,
	it satisfies for some $c_b>0$ and $s \in  {(0,1/2)}$
		\begin{equation} \label{eq:angularsing}
			\forall \, \theta \in (0, \pi/2], \quad \frac{c_b}{\theta^{1+2s}} \leq \sin \theta \, b(\cos \theta) \leq \frac{1}{c_b \,  \theta^{1+2s}};
		\end{equation}
	\item the kinetic factor $\Phi$ satisfies
		\begin{equation} \label{eq:Phi}
			\Phi(|v-v_*|)= |v-v_*|^\gamma \quad \text{with} \quad {\gamma>0},
		\end{equation}
		this assumption could be relaxed to assuming only that $\Phi$ satisfies $\Phi(\cdot) = C_\Phi \, |\cdot|^\gamma$ for some $C_\Phi>0$.
	\end{itemize}
We will consider $\mu$ the only global Maxwellian equilibrium of the equation with mass $1$, vanishing momentum and energy $3$:
	$$
	\mu(v):=(2 \pi)^{-3/2} e^{-|v|^2/2}.
	$$
We are interesting in the linearized operator around the equilibrium $\mu$ (not the whole nonlinear Boltzmann operator) which is defined at first order through
	$$
	\Lambda f := Q(\mu,f) + Q(f,\mu) - v \cdot \nabla_x f
	$$
and we thus consider the evolution equation~\eqref{FKE} with $L$ given by
\begin{equation} \label{eq:BO}
Lf := - (Q(\mu,f) + Q(f,\mu))
\end{equation}
with the collision operator $Q$ defined through~\eqref{eq:Q} and satisfying the conditions~\eqref{eq:B},~\eqref{eq:angularsing},~\eqref{eq:Phi}.

\smallskip
\subsection{Notations} We will denote $\langle w \rangle := (1+|w|^2)^{1/2}$ for any $w \in \R^d$. For convenience, we introduce
the following strictly positive operators
$$
\Lambda_v^2 :=   1 -\Delta_v,\qquad \Lambda_x^2 =  1 -\Delta_x
$$
 and the associated family of Fourier multipliers
 $$
\Lambda_x^\alpha := ( 1 -\Delta_x)^{\alpha/2}, \qquad  \Lambda_v^\beta :=  ( 1 -\Delta_v)^{\beta/2}, \qquad \qquad  \alpha, \beta \in \R
$$
 which act on a function in $\sss(\T^d\times \R^d)$ in the following way
 $$
 \widehat{ \Lambda_x^\alpha f} (\xi, \eta) = \sxi^\alpha \widehat{f}(\xi, \eta), \qquad \widehat{ \Lambda_v^\beta f} (\xi, \eta) = \seta^\beta \widehat{f}(\xi, \eta)
 $$
 where the hat corresponds to the Fourier transform in both $x$ (with corresponding variable $\xi \in \Z^3$) and $v$ (with corresponding variable $\eta \in \R^3$).  We also introduce the corresponding Sobolev spaces
 $$
 H^{\alpha,\beta}_{x,v} = \set{ f \in \sss', \, \, \Lambda_x^\alpha \Lambda_v^\beta f \in L^2},
 $$
 and we denote by $\norm{\cdot}_{\alpha,\beta}$ the corresponding norm defined by 
 $$
 \|f\|^2_{H^{\alpha,\beta}_{x,v}} := \sum_{\xi \in \Z^d} \int_{\eta \in \R^d} \langle \xi \rangle^{2\alpha} \langle \eta \rangle^{2 \beta} \widehat{f} (\xi,\eta)^2 \, d\eta. 
 $$ 
 Similarly, we introduce the weighted Sobolev spaces $H^{\alpha,\beta}_{x,v}(m)$ for $m$ a weight function (the typical example in the sequel will be $m(v) = \langle v \rangle^k$ for some $k \ge0$):
  $$
 H^{\alpha,\beta}_{x,v} (m)= \set{ f \in \sss', \, \, \Lambda_x^\alpha \Lambda_v^\beta (fm) \in L^2},
 $$
and we denote by $\norm{\cdot}_{H^{\alpha,\beta}_{x,v} (m)}$ the corresponding norm defined by
$$
 \|f\|^2_{H^{\alpha,\beta}_{x,v}(m)} := \|fm\|^2_ {H^{\alpha,\beta}_{x,v}}.
$$

We also define the classical weighted Sobolev space $ H^{n}_{x,v} (m)$, for $n\in \N$ by:

	$$
	\| f \|^2_{ H^{n}_{x,v}(m)} := \sum_{ |\alpha| \leq \ell, \,  |\beta| \leq n, \, |\alpha| + |\beta| \leq n}
	\| \partial^\alpha_v \partial^\beta_x (f m) \|^2_{L^2_{x,v}}.
	$$
{We use Fourier transform to define the general space $H^r_{x,v}(m)$ for $r \in \R^+$:
	\begin{equation} \label{homonorm}
	\|f\|^2_{H^r_{x,v}(m)} := \|f m\|^2_{H^r_{x,v}} = \sum_{\xi \in \Z^d} \int_{\eta \in \R^d} (1+|\xi|^2+|\eta|^2)^r \,|\widehat{fm}(\xi,\eta)|^2 \, d\eta
	\end{equation}
 where the hat still corresponds to the Fourier transform in $x$ and $v$.
In the case $r \in \N$, the norms given by the two previous formula are equivalent. We won't make any difference in the notation and will use one norm or the other at our convenience. It won't have any impact on our estimates since it will only add multiplicative universal constants.}

Let us remark that by classical results of interpolation (see for example~\cite{BookBL}), for any $r \in \R^+$, one can write
$$
H^r_{x,v}(m) = \left[H^{\lfloor r \rfloor}_{x,v}(m), H^{\lfloor r \rfloor+1}_{x,v}(m)\right]_{r-\lfloor r \rfloor,2}.
$$
The notation used above is the classical one of real interpolation. For sake of completeness, we briefly recall the meaning of this notation. For $C$ and $D$ two Banach spaces which are both embedded in the same topological separating vector space, for any $z \in A+B$, we define the $K$-function by
$$
K(t,z) := \inf_{z=c+d} \left(\|c\|_C+t\|d\|_D\right), \quad \forall \, t>0.
$$
We then give the definition of the space $[C,D]_{\theta,p}$ for $\theta \in (0,1)$ and $p \in [1,+\infty]$:
$$
[C,D]_{\theta,p} := \left\{z \in C+D, \, \, t \mapsto K(t,z)/t^\theta \in L^p\left(dt/t^{1/p}\right)\right\}.
$$

\Bk

\subsection{Main results and known results}
\subsubsection{Fractional Kolmogorov equation} With the notations introduced above, the fractional Kolmogorov equation reads
 $$
 \D_{t}f+v\cdot\nabla_{x}f + \Lambda_v^{2s}f =0
 $$
 and a natural question is wether $f$ benefits from some regularization
 induced by the  elliptic  properties of $\Lambda_v^{2s}$.
 The main result concerning the fractional Kolmogorov equation is the
 following:

 \begin{thm} \label{FKEthm}
 Let $r \in \R$ and $f$ be a solution of \eqref{FKE} with $L$ given by~\eqref{FKO} with initial data {$f_0 \in H^{r,0}_{x,v}$}. Then, there exists a constant $C_{r} >0$ independent of $f_0$ such that for all $t \in (0,1]$,  we have
$$
{ \norm{f(t)}_{r,s} \leq \frac{C_r}{t^{1/2}} \norm{f_0}_{r,0}}
\quad \text{and} \quad
 \norm{f(t)}_{r+s,0} \leq \frac{C_r}{t^{1/2+s}} \norm{f_0}_{r,0}.
 $$
 \end{thm}

 This result has already been proved in the case $s=1$ in \cite{Herau-JFA} by the first author and we give here a result concerning the cases $s \in (0,1)$ following essentially the same core of the method as there: we introduce a functional which is a Lyapunov functional for our equation for small times (see also~\cite{Vill-Hypo} by Villani and the references therein for this type of methods). From this property, we are then able to recover some regularization estimates quantified in time. Let us emphasize that the main difficulty for the fractional case is to find a good entropy function (it is of course not the same as in the non fractional case $s=1$). This type of result is of great use in the proof of the return to equilibrium in large functional spaces of solutions of inhomogeneous kinetic equations as in~\cite{MM-ARMA} by Mischler and Mouhot in the Fokker-Planck case ($s=1$) following a general method also presented in~\cite{MM-ARMA}. The homogeneous fractional Fokker-Planck case has been studied by the third author in \cite{Trist-CMS} where regularization properties in velocity (concerning the integrability of the solution) are investigated thanks to a fractional Nash inequality. 
 In Section~\ref{sec:FK}, devoted to the fractional Kolmogorov equation, we pay attention to give a proof for Theorem~\ref{FKEthm} \it without \rm using
 any kind of pseudodifferential tool (only Fourier multiplier).

 \subsubsection{Linearized Boltzmann equation without cutoff}
 In Section~\ref{sec:Bol}, what we aim to do is to prove some similar regularization properties for the linearized inhomogeneous Boltzmann equation without cutoff as we do for the fractional Kolmogorov equation (Theorem~\ref{FKEthm}).
We recall that the Boltzmann equation is of type~\eqref{FKE} in dimension $3$ with $L$ given by~\eqref{eq:BO}.
Here is the main result that we obtain on this model:
\begin{thm} \label{thm:LIBthm}
 Let { $r \in \N$,  $k' \ge 0$, $k >\max(\gamma/2 + 3+2\max(1,r)s, k'+\gamma+5/2)$} and $f$ be a solution of \eqref{FKE} with $L$ given by~\eqref{eq:BO} with initial data $f_0$. Then, there exists a constant $C_r>0$ independent of~$f_0$ such that we have the following regularization estimates. If { $f_0 \in H^{r,0}_{x,v}(\langle v \rangle^k)$ or $f_0 \in (H^{r,s}_{x,v}(\langle v \rangle^{k})'$}  where $(H^{r,s}_{x,v}(\langle v \rangle^{k}))'$ is the dual space of $H^{r,s}_{x,v}(\langle v \rangle^{k})$ with respect to $H^{r,0}_{x,v}(\langle v \rangle^{k})$, for any $t \in (0,1]$,
 $$
{ \|f(t)\|_{H^{r,s}_{x,v}(\langle v \rangle^{k'})} \le \frac{C_r}{t^{1/2}} \|f_0\|_{H^{r,0}_{x,v}(\langle v \rangle^{k})} \quad \text{or} \quad \|f(t)\|_{H^{r,0}_{x,v}(\langle v \rangle^{k'})} \leq \frac{C_r}{t^{1/2}} \|f_0\|_{(H^{r,s}_{x,v}(\langle v \rangle^{k}))'}}.
 $$
If $f_0 \in H^{r,0}_{x,v}(\langle v \rangle^k)$ or $f_0 \in (H^{r+s,0}_{x,v}(\langle v \rangle^{k}))'$, where $(H^{r+s,0}_{x,v}(\langle v \rangle^{k}))'$ is the dual space of $H^{r+s,0}_{x,v}(\langle v \rangle^{k})$ with respect to $H^{r,0}_{x,v}(\langle v \rangle^{k})$, we also have for any $t \in (0,1]$,
$$
{\|f(t)\|_{H^{r+s,0}_{x,v}(\langle v \rangle^{k'})} \le \frac{C_r}{t^{1/2+s}} \|f_0\|_{H^{r,0}_{x,v}(\langle v \rangle^{k})}
\quad \text{or} \quad \|f(t)\|_{H^{r,0}_{x,v}(\langle v \rangle^{k'})} \leq \frac{C_r}{t^{1/2+s}} \|f_0\|_{(H^{r+s,0}_{x,v}(\langle v \rangle^{k}))'}}.
 $$

\end{thm}
%

First, we have to underline that it is the first result of regularization quantified in time on the Boltzmann equation without cutoff and that it is a key point for the development of the Cauchy theory of perturbative solutions in \cite{HTT1*} by the same three authors for the nonlinear Boltzmann equation without cutoff (the condition on the power $k$ comes from this paper).

The singularity of the Boltzmann kernel in the non cutoff case implies that the Boltzmann operator without cutoff (that we will describe later on) behaves as a fractional Laplacian in velocity:
$$
Q(g,f) \approx - C_g (-\Delta_v)^s f + \, \text{lower order terms}
$$
with $C_g$ depending only on the physical properties of $g$.
 This type of result has already been studied in the homogeneous and non-homogeneous cases. The gain in velocity is quite obvious to observe even if it is complicated to understand it precisely: up to now, the most common way to understand it is through an anisotropic norm (see~\cite{GS-JAMS} by Gressman and Strain and~\cite{AMUXY-ARMA2} by Alexandre et al.). It is then natural to expect that the transport term allows to transfer the gain in velocity to the space variable. We refer to the references quoted in \cite{AHL*} for a review of this type of properties. Let us mention that the paper~\cite{AHL*} by Alexandre et al. is the first one in which the hypoellipticity features of the operator have been deeply analyzed.

Our strategy here is to use the same method as for Kolmogorov type equations introduced in~\cite{Herau-JFA} by the first author. In short, except from the fact that the use of pseudodifferential tools is required and thus there are many additional technical difficulties, the method is the same as for the fractional Kolmogorov equation.
For purposes of comparison, we can also mention that this kind of strategy has also been applied successfully to the Landau equation in \cite{CTW-ARMA} by Carrapatoso et al..
However, the study of this kind of properties is much harder in the case of the Boltzmann equation without cutoff since the gain in regularity is less clear and consists in an anisotropic gain of fractional derivatives: we have to exploit the fact that one can write a part the Boltzmann linearized operator as a pseudo-differential operator, in the spirit of what has been done in~\cite{AHL*}.

 Indeed, we adapt here some ideas from there allowing to do computations for operators - including the Boltzmann one - whose symbols are in an adapted class called here $S_K$, where  $K$ is a large parameter. Let us point out that those classes are complicated partly because the order of the symbols does not decrease with derivation, which induces some great technical difficulties.
The computations are done using the Wick quantization, widely studied in particular by Lerner (see \cite{LernerBook1} and \cite{LernerBook2}), which has very nice positivity properties. This allows to adapt  to the Boltzmann case the Lyapunov strategy already introduced in~\cite{Herau-JFA} for the Kolmogorov case and in the second section of this article for the fractional Fokker-Planck case.

 It is also important to underline the fact that this pseudo-differential study is not done on the whole linearized operator but only on a well-chosen part of it (this is the object of Subsection~\ref{subsec:reg}). Indeed, thanks to Duhamel formula, we will then be able to recover an estimate on the whole semigroup, the one associated to $\Lambda$ (see Lemma~\ref{lem:lambda_1/lambda}). \bigskip

\noindent\textbf{Acknowledgments.} This research has been supported by the \'Ecole Normale Sup\'erieure through the project {\it Analyse de solutions d'\'equations de la th\'eorie cin\'etique des gaz}. The first author thanks the Centre Henri Lebesgue ANR-11-LABX-0020-01 for its support and the third author thanks the ANR EFI:  ANR-17-CE40-0030.

\bigskip
\subsection{Outline of the paper}
In Section~\ref{sec:FK}, we prove Theorem~\ref{FKEthm} on the fractional Kolmogorov equation. In Section~\ref{sec:Bol}, we focus on the linearized Boltzmann equation and more precisely on Theorem~\ref{thm:LIBthm} and most of its proof. The remainder of which contains pseudo-differential arguments and will be found in Section~\ref{sec:quantization}. In Section~\ref{sec:general}, we explain how to generalize the result for higher order derivatives in the Boltzmann case. Finally, in the Appendix, we present the pseudo-differential tools that we shall use in full generality.

\bigskip
\section{The fractional Kolmogorov equation} \label{sec:FK}

This section is devoted to the proof of Theorem \ref{FKEthm}.
Following \cite{Herau-JFA}, we shall show below that this entropy functional is decreasing with time, and this will imply the result in the last subsection.

\subsection{A Lyapunov functional}  \label{sublyapKFP}

Let $f$ be a solution of \eqref{FKE} with $L$ given by \eqref{FKO} and with initial data $f_0$. We first deal with the case $r=0$. We follow the lines of the proof given in~\cite{Herau-JFA} and introduce an adapted entropy
functional defined  for all $t\geq 0$ by
$$
\hhh(t) := C\norm{f}^2 + D t \norm{\Lambda_v^{s-1} \nabla_v f}^2 + E t^{1+s} \Re \sep{\Lambda_v^{s-1} \nabla_v f, \Lambda_x^{s-1} \nabla_x f} + t^{1+2s} \norm{ \Lambda_x^{s-1} \nabla_x f}^2
$$
for large constants $C$, $D$, $E$ to be chosen later, where $\norm{\cdot}$ is the usual $L^2$ norm, $\sep{ \cdot,\cdot}$ is the usual (complex) $L^2$ scalar product, we also have denoted
$$
\norm{\Lambda_w^{s-1} \nabla_w f}^2 = \sum_{k=1}^d \norm{\Lambda_w^{s-1} \partial_{w_k} f}^2 \quad \text{for $w=x,v$}
$$
and
$$
\sep{\Lambda_v^{s-1} \nabla_v f, \Lambda_x^{s-1} \nabla_x f} = \sum_{k=1}^d \sep{\Lambda_v^{s-1} \partial_{v_k} f, \Lambda_x^{s-1} \partial_{x_k} f}.
$$

The first step in the study is to show
that $\hhh$ is indeed non-negative. The lemma below shows in addition that for all $t\ge 0$,
$\hhh(t)$ controls the $H^s$ norm (where $H^{s} := H^{0,s} \cap H^{s,0}$).

\begin{lem} \label{hh1}
If $E \leq \sqrt{D}$ then for all  $t\geq 0$ we have $\hhh(t) \geq 0$. Precisely we have
\begin{equation*}
0 \leq  C\norm{f}^2 + \frac{D}{2} t \norm{\Lambda_v^{s-1}\nabla_v f}^2  + \frac{1}{2} t^{1+2s} \norm{ \Lambda_x^{s-1} \nabla_x f}^2
\leq \hhh(t).
\end{equation*}
\end{lem}
\preuve The proof is direct using the time-dependant Cauchy-Schwarz inequality
$$
E t^{s} \abs{ \sep{\Lambda_v^{s-1} \nabla_v f, \Lambda_x^{s-1} \nabla_x f}}  \leq
\frac{E^2}{2}  \norm{\Lambda_v^{s-1} \nabla_v f}^2  + \frac{1}{2} t^{2s} \norm{ \Lambda_x^{s-1} \nabla_x f}^2.
$$
\fin

The main ingredient in the proof of Theorem \ref{FKEthm} is the following commutation equality: for $j \in \llbracket 1,d \rrbracket$,
$$
\adf{ \D_{v_j}, v_j \D_{x_j}} = \D_{x_j}.
$$
In the same spirit, we shall need later the following lemma
giving formulas for  slightly modified commutators. We denote from now on $X := v \cdot \nabla_x$ the Vlasov operator and $X_j := v_j \D_{x_j}$, so that $X = \sum_{j=1}^d X_j$ and the previous fundamental equality reads $ \adf{ \D_{v_j}, X_j} = \D_{x_j}$.

\begin{lem} \label{commm} For $k \in \llbracket 1,d \rrbracket$, we have
$$
\adf{ \Lambda_v^{s-1} \D_{v_k}, X} =  \Lambda_v^{s-1}  \D_{x_k} +(1-s) \D_{v_k} \sum_{j=1}^d \D_{v_j}  \Lambda_v^{s-3}  \D_{x_j}
$$
  and 
$$
\adf{ \Lambda_v^{s-1} \D_{v_k}, \Lambda_v^{2s}} = \adf{ \Lambda_x^{s-1} \D_{x_k}, \Lambda_v^{2s}} = \adf{ \Lambda_x^{s-1} \D_{x_k}, X} = 0.$$
\end{lem}

\preuve
For the three last equalities, the result is immediate since differentiation in velocity and spatial direction commute. Let us deal with the first one. Let $j$, $k \in \llbracket 1,d \rrbracket$. We check that
that the commutator $\adf{ \Lambda_v^{s-1} \D_{v_k}, X_j}$ is in fact a Fourier multiplier whose symbol reads
\begin{equation*} 
\begin{split}
\sigma \sep{ \adf{ \Lambda_v^{s-1} \D_{v_k}, X_j}} & = \frac{1}{i} \set{ \seta^{s-1} i\eta_k, iv_j \xi_j}
\end{split}
\end{equation*}
where we denote by $\set{ \cdot, \cdot}$ the Poisson bracket of two functions. Let us mention that in the Fourier formalism, we have that for $\alpha \in \R$,
\begin{eqnarray*}
\sigma( \D_{v_k}) = i\eta_k, \quad \sigma(\D_{x_j}) = i\xi_j, \qquad \sigma(X_j) = iv_j\xi_j, \quad \sigma(\Lambda_v^\alpha) = \seta^\alpha,\quad \sigma( -\Delta_v) = |\eta|^2.
\end{eqnarray*}
We then have 
$$
\sigma  \sep{ \adf{ \Lambda_v^{s-1} \D_{v_k}, X_j}}  = i \xi_j \set{ \seta^{s-1} \eta_k, v_j}
 =  i\xi_j \sep{ (s-1) \eta_j \eta_k \seta^{s-3} + \seta^{s-1} \delta_{kj}}
$$
where $\delta_{kj}$ is the Kronecker delta of $(k,j)$. Then, summing on $j$, we get:
\begin{align*}
\sigma  \sep{ \adf{ \Lambda_v^{s-1} \D_{v_k}, X}} &= i \sum_{j=1}^d \xi_j \sep{ (s-1) \eta_j \eta_k \seta^{s-3} + \seta^{s-1} \delta_{kj}} \\
&= i  \seta^{s-3} \left(\xi_k (1+s\eta_k^2) + \sum_{j \neq k} \eta_j ((s-1)\xi_j \eta_k + \xi_k \eta_j) \right).
\end{align*}
Coming back on the non-Fourier side, we obtain:
\begin{align*}
\adf{ \Lambda_v^{s-1} \D_{v_k}, X} &=(1-s\D_{v_k}^2)\Lambda_v^{s-3}  \D_{x_k} +(1-s) \D_{v_k} \sum_{j \neq k} \D_{v_j}  \Lambda_v^{s-3}  \D_{x_j} -   \sum_{j \neq k} \D_{v_j}^2 \Lambda_v^{s-3}  \D_{x_k} \\
&=(1 - \Delta_v) \Lambda_v^{s-3}  \D_{x_k} +(1-s) \D_{v_k}^2 \Lambda_v^{s-3}  \D_{x_k} +(1-s) \D_{v_k} \sum_{j \neq k} \D_{v_j}  \Lambda_v^{s-3}  \D_{x_j} \\
&=  \Lambda_v^{s-1}  \D_{x_k} +(1-s) \D_{v_k} \sum_{j=1}^d \D_{v_j}  \Lambda_v^{s-3}  \D_{x_j}
\end{align*}
which is the required result.
\fin

We now show that $\hhh$ is indeed a Lyapunov function (entropy functional).
\begin{lem} \label{deriv}  For well chosen (arbitrarily large) constants $C$, $D$ and $E$ we have
$$
\ddt \hhh(t) \leq 0, \quad \forall \, t \in [0,1].
$$
\end{lem}

\preuve Using the previous lemma, we shall compute the time derivative of each terms appearing in the definition of $\hhh$. For convenience  we introduce the operator associated the the Kolmogorov equation
$$
P := X + \Lambda_v^{2s}
$$
so that $f$ satisfies $\D_t f+ Pf = 0 $. We do below all the computations in (the complex) $L^2$.

We first notice that
 \begin{equation*} 
 \ddt \norm{f}^2 = -2\Re (Pf,f) = -2\Re ( (X + \Lambda^{2s}_v)f,f) =  -2  (\Lambda_v^{2s} f,f)
 \end{equation*}
 since $X$ is skew-adjoint. Using Parseval formula on the right-hand side we get that the first term in the derivative of $\hhh$ is
  \begin{equation} \label{f}
 \ddt C\norm{f}^2 =  -\big(\underbrace{2 C \seta^{2s}}_{I} \hf, \hf \big).
 \end{equation}
 Note that this term is non-positive.

For the second term in the derivative of $\hhh$, we have
\begin{equation*}
\ddt \sep{  t\norm{\Lambda_v^{s-1} \nabla_v f}^2 } =  \sum_{k=1}^d \left(\norm{\Lambda^{s-1}_v \D_{v_k} f}^2 + t \ddt (\Lambda^{s-1}_v \D_{v_k} f , \Lambda^{s-1}_v \D_{v_k} f) \right).
 \end{equation*}
Let us fix $k \in \llbracket 1, d \rrbracket$. The derivative of the $k$-th term in the last term writes
\begin{align*}
\ddt \norm{\Lambda^{s-1}_v \D_{v_k} f}^2  = & -2 \Re (\Lambda_v^{s-1} \D_{v_k} P f , \Lambda_v^{s-1} \D_{v_k} f)  \\
  = &  -2 \Re ( \Lambda_v^{s-1} \D_{v_k} \Lambda_v^{2s} f , \Lambda_v^{s-1} \D_{v_k} f) -2 \Re (\Lambda_v^{s-1} \D_{v_k}  X f , \Lambda_v^{s-1} \D_{v_k} f) \\
   = & -2 \Re (\Lambda_v^{2s} \Lambda_v^{s-1} \D_{v_k}  f , \Lambda_v^{s-1} \D_{v_k} f) -2 \Re (X \Lambda_v^{s-1} \D_{v_k}   f , \Lambda_v^{s-1} \D_{v_k} f) \\
  & -2 \Re (\adf{ \Lambda_v^{s-1} \D_{v_k} , \Lambda_v^{2s}} f , \Lambda_v^{s-1} \D_{v_k} f) -2 \Re (\adf{ \Lambda_v^{s-1} \D_{v_k}, X  }  f , \Lambda_v^{s-1} \D_{v_k} f) \\
  = & -2 \Re (\Lambda_v^{2s} \Lambda_v^{s-1} \D_{v_k}  f , \Lambda_v^{s-1} \D_{v_k} f)  \\[-3mm]
  &\quad - 2 \Re ( \Lambda_v^{s-1}  \D_{x_k} f+(1-s) \D_{v_k} \sum_{j=1}^d \D_{v_j}  \Lambda_v^{s-3}  \D_{x_j} f, \Lambda_v^{s-1} \D_{v_k} f) \\
   = &  2\Re (\Lambda_v^{4s-2} \D_{v_k}^2  f ,  f)  + 2 \Re (  \Lambda_v^{2s-2}  \D_{x_k} \D_{v_k} f  , f) \\[-2mm]
  &\quad +2\Re ((1-s) \D^2_{v_k} \sum_{j=1}^d \D_{v_j}  \Lambda_v^{2s-4}  \D_{x_j} f , f)
 \end{align*}
where we used that $X$ is skew-adjoint and the commutation expressions in  Lemma \ref{commm}. Writing the right-hand side on the Fourier side, summing over $k$ and using Cauchy-Schwarz inequality gives us:
  \begin{align*}
  \ddt \norm{\Lambda^{s-1}_v \nabla_v f}^2
  & \leq - 2 (\seta^{4s-2} |\eta|^2 \hf, \hf) + 2(2-s) (\seta^{2s-1} \sxi \hf, \hf) \\
  & \leq - 2(\seta^{4s}  \hf, \hf) +2(\seta^{4s-2}  \hf, \hf) + 2(2-s) ( \seta^{2s-1}  \sxi \hf, \hf).
  \end{align*}
 The second term in $\hhh$ therefore satisfies
\begin{multline} \label{vvf}
 \ddt \sep{ D t \norm{\Lambda^{s-1}_v \nabla_v f}^2 } \\ \leq \bigg( \big( \underbrace{D \seta^{2s}}_{i} - \underbrace{2Dt\seta^{4s}}_{II} + \underbrace{2Dt\seta^{4s-2}}_{ii} + \underbrace{2(2-s) Dt  \seta^{2s-1} \sxi}_{iii}\big) \hf, \hf \bigg).
 \end{multline}
 We note that the term corresponding to II is non-positive, and that the three other ones are non-negative.
We now deal with the third term in the derivative of $\hhh$:

 \begin{align*}
&\quad \ddt \bigg(  t^{1+s} \Re\sep{\Lambda_v^{s-1} \D_v f, \Lambda_x^{s-1} \D_x f} \bigg) \\
&= (1+s) t^{s} \Re \sep{\Lambda_v^{s-1} \nabla_v f, \Lambda_x^{s-1}\nabla_x f} + t^{1+s} \ddt  \Re \sep{\Lambda_v^{s-1} \nabla_v f, \Lambda_x^{s-1} \nabla_x f} \\[-0.5mm]
&=  (1+s) t^{s} \Re \sep{\Lambda_v^{s-1} \nabla_v f, \Lambda_x^{s-1}\nabla_x f} + t^{1+s} \sum_{k=1}^d\ddt  \Re \sep{\Lambda_v^{s-1} \D_{v_k} f, \Lambda_x^{s-1} \D_{x_k} f}.
 \end{align*}
The $k$-th derivative in the last term writes
\begin{equation*}\begin{split}
&\quad \ddt \Re(\Lambda_v^{s-1} \D_{v_k} f , \Lambda_x^{s-1} \D_{x_k} f)  \\
=  &- \Re (\Lambda_v^{s-1} \D_{v_k} P f , \Lambda_x^{s-1} \D_{x_k} f)- \Re (\Lambda_v^{s-1} \D_{v_k}  f , \Lambda_x^{s-1} \D_{x_k} P f)  \\
  = &  -2 \Re (\Lambda_v^{2s} \Lambda_v^{s-1} \D_{v_k}  f , \Lambda_x^{s-1} \D_{x_k} f) \\
    & - \Re (\Lambda_v^{s-1} \D_{v_k}  X f , \Lambda_x^{s-1} \D_{x_k} f)- \Re (\Lambda_v^{s-1} \D_{v_k}  f , \Lambda_x^{s-1} \D_{x_k} X f) \\
    = &  -2 \Re (\Lambda_v^{2s} \Lambda_v^{s-1} \D_{v_k}  f , \Lambda_x^{s-1} \D_{x_k} f) \\
    & - \Re (\adf{ \Lambda_v^{s-1} \D_{v_k} , X} f , \Lambda_x^{s-1} \D_{x_k} f)- \Re (\Lambda_v^{s-1} \D_{v_k}  f , \adf{ \Lambda_x^{s-1} \D_{x_k}, X } f) \\
    & - \Re (X \Lambda_v^{s-1} \D_{v_k}   f , \Lambda_x^{s-1} \D_{x_k} f)- \Re (\Lambda_v^{s-1} \D_{v_k}  f , X \Lambda_x^{s-1} \D_{x_k}  f).
    \end{split}
    \end{equation*}
   We use again that $X$ is skew-adjoint and observe that it implies that the sum of the last two terms is zero by compensation. The previous term is also zero since the commutator inside is zero. With Lemma \ref{commm}, we obtain
    \begin{align*}
    \ddt \Re (\Lambda_v^{s-1} \D_{v_k} f , \Lambda_x^{s-1} \D_{x_k} f) = & -2 \Re (\Lambda_v^{2s} \Lambda_v^{s-1} \D_{v_k}  f , \Lambda_x^{s-1} \D_{x_k} f) \\
    &- \Re \big(  \Lambda_v^{s-1}  \D_{x_k} f+(1-s) \D_{v_k} \sum_{j=1}^d \D_{v_j}  \Lambda_v^{s-3}  \D_{x_j}  f , \Lambda_x^{s-1} \D_{x_k} f\big).
 \end{align*}
  Writing the right-hand side on the Fourier side then gives
 \begin{equation*}\begin{split}
\ddt \Re(\Lambda_v^{s-1} \D_{v_k} f , \Lambda_x^{s-1} \D_{x_k} f)
  & = -  (\seta^{3s-1} \sxi^{s-1} \eta_k \xi_k \hf, \hf) -  (\seta^{s-1}  \xi_k^2\sxi^{s-1}  \hf, \hf) \\
  &\quad + (1-s) \sum_{j=1}^d(\eta_k \seta^{s-3} \xi_k \sxi^{s-1}  \eta_j \xi_j \hf, \hf ).
 \end{split}
 \end{equation*}
 Then, taking the sum overs $k$ gives us (using Cauchy-Schwarz inequality):
 \begin{align*}
 \ddt \Re(\Lambda_v^{s-1} \nabla_v f , \Lambda_x^{s-1}\nabla_x f)
 &\le d (\seta^{3s} \sxi^{s}   \hf, \hf) - ( \seta^{s-1}  |\xi|^2\sxi^{s-1}  \hf, \hf) \\
  &\quad+ (1-s) \bigg( \bigg(\sum_{j=1}^d  \eta_j \xi_j \bigg)^2 \seta^{s-3} \sxi^{s-1}\hf, \hf \bigg) \\
&\le d (\seta^{3s} \sxi^{s}   \hf, \hf) - (\seta^{s-1}  |\xi|^2\sxi^{s-1}  \hf, \hf) \\
 &\quad + (1-s) ( |\eta|^2 |\xi|^2 \seta^{s-3} \sxi^{s-1}\hf, \hf ) \\
& \le d (\seta^{3s} \sxi^{s}   \hf, \hf) - s (\seta^{s-1} \sxi^{s+1} \hf, \hf) +  (\seta^{s-1} \sxi^{s-1} \hf, \hf).
 \end{align*}
 We therefore get that the third term in $\hhh$ satisfies:
\begin{multline} \label{vxf}
 \ddt \bigg( E  t^{1+s} \Re\sep{\Lambda_v^{s-1} \nabla_v f, \Lambda_x^{s-1} \nabla_x f} \bigg) \\ \leq \bigg( \big( \underbrace{E(s+1)t^s \seta^{s} \sxi^{s} }_{iv}   + \underbrace{d E  t^{1+s}\seta^{3s}\sxi^{s}}_{v} \\ - \underbrace{E s t^{1+s} \seta^{s-1}\sxi^{s+1}}_{III} + \underbrace{E t^{1+s}  \seta^{s-1}  \sxi^{s-1}}_{vi}    \big) \hf, \hf \bigg).
 \end{multline}
 We note that the term corresponding to III is non-positive, and that the three other ones are non-negative.

We can now deal with the last term in the derivative of $\hhh$.
We write

\begin{align*}
\ddt \sep{  t^{1+2s}\norm{\Lambda_x^{s-1} \nabla_x f}^2 } &=  (1+2s) t^{2s} \norm{\Lambda_x^{s-1} \nabla_x f}^2 + t^{1+2s} \ddt \norm{\Lambda_x^{s-1} \nabla_x f }^2 \\
&= (1+2s) t^{2s} \sum_{k=1}^d \norm{\Lambda_x^{s-1} \D_{x_k} f}^2 +t^{1+2s} \sum_{k=1}^d \ddt \norm{\Lambda_x^{s-1} \D_{x_k} f }^2.
 \end{align*}
The $k$-th derivative of the last term writes
\begin{equation*}\begin{split}
\ddt \norm{\Lambda_x^{s-1} \D_{x_k} f }^2  = & -2 \Re (\Lambda_x^{s-1} \D_{x_k} P f , \Lambda_x^{s-1} \D_{x_k} f)  \\
  = &  -2 \Re (\Lambda_v^{2s} \Lambda_x^{s-1} \D_{x_k}  f , \Lambda_x^{s-1} \D_{x_k} f) -2 \Re (X \Lambda_x^{s-1} \D_{x_k}  f , \Lambda_x^{s-1} \D_{x_k} f) \\
  = & -2 \Re (\Lambda_v^{2s} \Lambda_x^{s-1} \D_{x_k}  f , \Lambda_x^{s-1} \D_{x_k} f).  \\
 \end{split}
 \end{equation*}
 We used here the last commutations properties in Lemma \ref{commm} and again that $X$ is skew-adjoint. Writing the right-hand side on the Fourier side  and summing on $k$ gives
 \begin{align*}
\ddt \norm{\Lambda_x^{s-1} \nabla_x f }^2
  & = - 2(\seta^{2s} \sxi^{2s-2} |\xi|^2 \hf, \hf) \\
  & = - 2(\seta^{2s} \sxi^{2s}  \hf, \hf) +2(\seta^{2s} \sxi^{2s-2} \hf, \hf).
 \end{align*}
 The fourth term in the derivative of  $\hhh$ therefore satisfies:
\begin{multline} \label{xxf}
 \ddt \sep{ t^{1+2s} \norm{\Lambda_x^{s-1} \nabla_x f }^2 }  \\
 \leq \bigg( \big( \underbrace{(1+2s) t^{2s} \sxi^{2s}}_{vii} - \underbrace{2t^{1+2s} \seta^{2s}\sxi^{2s}}_{IV} + \underbrace{2t^{1+2s} \seta^{2s}\sxi^{2s-2}}_{viii}\big) \hf, \hf \bigg).
 \end{multline}
 We note that the term corresponding to IV is non-positive, and that the other ones are non-negative.

\bigskip
Now we look at the different
terms appearing in formulas (\ref{f}-\ref{xxf}). We want to show that with a good choice of constants $C$, $D$ and $E$, the corresponding sum is non-positive, and therefore $\hhh$ is indeed a Lyapunov functional. We shall study each non-negative term (small letters $(i)$ to $(viii)$ ) and show that they can be controlled by combinations of terms $I$ to $IV$, using essentially the Hölder inequality in $\R^2$. We restrict the study to $t \in [0,1]$.

The terms $(i)$  and $(ii)$ are immediately bounded by $I/10$ if
\begin{equation} \label{condi}
2D \leq 2C/10.
\end{equation}
 since $s\leq 1$.
The term $(iii)$ is a little bit trickier. We check that for any $\eps_{iii} >0$
$$
t \seta^{2s-1} \sxi \leq \eps_{iii}^{-1} \seta^{2s} + \eps_{iii}^s t^{1+s} \seta^{s-1}\sxi^{s+1}.
$$
Multiplying this inequality by $2D$ implies that $(iii) \leq I/10 + III/10$ if the following conditions are satisfied
\begin{equation} \label{condiii}
\eps_{iii}^{-1} 2(2-s)D \leq 2C/10,  \qquad \eps_{iii}^s 2(2-s)D \leq Es/10.
\end{equation}
We now deal with the term $(iv)$. We first check that for any $\eps_{iv} >0$
$$
t^s \seta^{s} \sxi^s \leq \eps_{iv}^{-1} \seta^{2s} + \eps_{iv}^{1/s} t^{1+s} \seta^{s-1}\sxi^{s+1}.
$$
Multiplying this inequality by $E(s+1)$ implies that $(iv) \leq I/10 + III/10$ if the following conditions are satisfied
\begin{equation} \label{condiv}
\eps_{iv}^{-1} E(s+1) \leq 2C/10,  \qquad \eps_{iv}^{1/s} E(s+1) \leq Es/10.
\end{equation}
For the term $(v)$, we also have to give a refined estimate.
We first check that for any $\eps_{v} >0$
$$
t^{1+s} \seta^{3s} \sxi^s \leq \eps_{v}^{-1} t \seta^{4s} + \eps_{v} t^{1+2s} \seta^{2s}\sxi^{2s}.
$$
Multiplying this inequality by $dE$ implies that $(v) \leq II/10 + IV/10$ if the following conditions are satisfied
\begin{equation} \label{condv}
\eps_{v}^{-1} dE \leq 2D/10,  \qquad \eps_{v} dE \leq 2/10.
\end{equation}
The term $(vi)$ is easily handled since $s\leq 1$, and we directly get that
$(vi) \leq I/10$ if
\begin{equation} \label{condvi}
E \leq 2C/10.
\end{equation}
Now we study the term $(vii)$. We first notice that for any $\eps_{vii} >0$
$$
t^{2s}  \sxi^{2s} \leq \eps_{vii}^{-1}  \seta^{2s} + \eps_{vii}^{\frac{1-s}{2s}} t^{1+s} \seta^{s-1}\sxi^{s+1}.
$$
Multiplying this inequality by $(1+2s)$ implies that $(vii) \leq I/10 + III/10$ if the following conditions are satisfied
\begin{equation} \label{condvii}
\eps_{vii}^{-1} (1+2s) \leq 2C/10,  \qquad \eps_{vii}^{\frac{1-s}{2s}} (1+2s) \leq Es/10.
\end{equation}
To finish, the term $(viii)$ is also easily handled since $s\leq 1$, and we directly get that
$(viii) \leq I/10$ if
\begin{equation} \label{condviii}
2 \leq 2C/10.
\end{equation}

Now we can do the synthesis and check that we can choose (in order  of reverse appearance) the constants $C$, $D$, $E$ and the small constants
$\eps_{iii}$, $\eps_{iv}$, $\eps_{v}$ and $\eps_{vii}$ such that
conditions (\ref{condi}-\ref{condviii}) are satisfied. Note that $D$ and after that $C$  can be taken arbitrarily large at the end of this procedure. We obtain therefore that
\begin{equation} \label{synthpartielle}
\ddt \hhh(t) \leq -\frac{1}{10} \sep{(I + II + III + IV)\hf, \hf} \leq 0
\end{equation}
and the proof is complete. \fin

%

Then we are able to conclude the proof of the main result Theorem~\ref{FKEthm} concerning  the fractional Kolmogorov equation.

\subsection{Proof of Theorem \ref{FKEthm}}
We first prove the result for $r=0$.
Let $C$, $D$ and $E$ be constants given by Lemmas~\ref{hh1} and~\ref{deriv} and let us take $f_0 \in \sss$.
From Lemma \ref{deriv}, we first get that for all $t \in [0,1]$
$$
\hhh(t) \leq \hhh(0) = C \norm{f_0}^2.
$$
Using now Lemma \ref{hh1}, we get in particular
$$
\frac{D}{2} t\norm{ \Lambda_v^s f}^2 \leq C\norm{f}^2 + \frac{D}{2} t \norm{\Lambda_v^{s-1} \D_v f}^2 \leq \hhh(t) \leq C\norm{f_0}^2
$$
and this implies the result for the velocity regularization.
Similarly, using again Lemma~\ref{hh1}, we have
$$
\frac{1}{2} t^{1+2s} \norm{ \Lambda_x^{s}  f}^2 \leq C\norm{f}^2 +\frac{1}{2}t^{1+2s} \norm{ \Lambda_x^{s-1} \D_x f}^2 \leq \hhh(t) \leq C\norm{f_0}^2
$$
and this gives the regularization result for $r=0$ in the spatial direction.

For $r\in \R$, we just use the fact that $P$ commutes with
$\Lambda_x^r$ which implies that for $f$  solution of $\D_t f + Pf = 0$ with initial data $f_0$,   $\Lambda_x^r f$ is the solution of
$\D_t \Lambda_x^r f + P \Lambda_x^r f = 0$ with initial data $\Lambda_x^r f_0$. We can therefore apply the result on $L^2$ to $\Lambda_x^r f$ and this directly gives that
$$
\frac{D}{2} t\norm{ \Lambda_v^s\Lambda_x^r f}^2  \leq C\norm{\Lambda_x^r f_0}^2 \textrm{ and } \frac{1}{2} t^{1+2s} \norm{ \Lambda_x^{s+r}  f}^2 \leq  C\norm{\Lambda_x^r f_0}^2.
$$
This provides us the estimates for any $r\in \R$ and $f_0 \in \sss$. The general result for initial data in the corresponding spaces follows by density of $\sss$. The proof is complete. \fin

\bigskip
\section{The inhomogeneous Boltzmann without cutoff case} \label{sec:Bol}
This section is devoted to the study of the Boltzmann equation without cutoff case and more precisely, to the proof of Theorem~\ref{thm:LIBthm} (we recall that in this section, the dimension $d$ equals $3$).
We start by making a few comments on this theorem:
\begin{itemize}

\item The result is not optimal in the sense that there is a loss in weight in our estimates. But we strongly believe that one could obtain a better estimate (concerning the weights) carrying out a more careful study of the operator~$\Lambda$. Indeed, in our proof, we perform a rough splitting of it and we  use Duhamel formula to recover an estimate on the whole semigroup $S_{\Lambda}(t)$. We could have not split the operator and study it completely, that would certainly provides us a better result. However, the proof would be much more complicated and we are here interested in the gain of regularity in terms of derivatives (not in terms of weights). Furthermore, the result that we obtain is enough to develop our perturbative Cauchy theory in~\cite{HTT1*} because we have some leeway in the weights in our proof.

\item Another important fact is that Theorem~\ref{thm:LIBthm} provides a {``primal'' and a ``dual'' result of regularization, roughly speaking, from $L^2$ into $H^s$ and from $H^{-s}$ into~$L^2$. The fact that we also develop a dual result is directly related to the use of this theorem that we make in~\cite{HTT1*}. We will only present the proof of the dual result into full details, we just explain how to adapt it in the primal case (which is easier to handle) in Section~\ref{sec:general}. }

\end{itemize}

\subsection{Splitting of the operator {for the dual result}} \label{subsec:splitting}
As already mentioned above, we are going to study the regularization properties only of a part of $\Lambda$, we thus start by exhibiting a splitting of it.
There are at least two types of splittings that one can consider to separate grazing and non-grazing collisions, depending on the adopted troncature function: one can cut the small $\theta$ or the small $|v'-v|$. For our purpose, we will work with the second option which is more adapted to the study of hypoelliptic properties of the Boltzmann collision operator. To do that we introduce the truncation function $\chi \in \mathcal{D}(\R)$ which satisfies $\mathds{1}_{[-1,1]} \le \chi \le \mathds{1}_{[-2,2]}$ and $\chi_{\delta} (\cdot) := \chi(\cdot/\delta)$ for $\delta>0$ and consider the troncature function $\chi_\delta(|v'-v|)$. We denote $Q_\delta$ the operator associated to the kernel:
$$
B_\delta (v-v_*,\sigma) := \chi_\delta(|v'-v|) \, b(\cos \theta) \, |v-v_*|^\gamma
$$
and $Q_\delta^c$ the one associated to the remainder part of the kernel:
$$
B^c_\delta (v-v_*,\sigma):=(1-\chi_\delta(|v'-v|)) \, b(\cos \theta) \, |v-v_*|^\gamma.
$$
We then have:
$$
\begin{aligned}
\Lambda f &= - v \cdot \nabla_x f + Q_\delta (\mu, f) + Q_\delta^c(\mu,f) + Q(f,\mu) \\
&= \Bigg(- K \langle v \rangle^{\gamma+2s} f- v \cdot \nabla_x f + \int_{\R^3 \times \mathbb{S}^2} B_\delta(v-v_*,\sigma) (\mu_* f' - \mu'_*f)\, \d\sigma \, {\rm d} v_* \Bigg)\\
&\quad+\Bigg(  K \langle v \rangle^{\gamma+2s} f + \int_{\R^3 \times \mathbb{S}^2} B_\delta(v-v_*,\sigma) (\mu'_* - \mu_*)(f'+f) \, \d\sigma \, \dv_*\\
&\hskip 4.5cm+ \int_{\R^3\times \mathbb{S}^2} B_\delta^c(v-v_*,\sigma) (\mu'_* f' - \mu_* f) \, \d\sigma \, \dv_* + Q(f,\mu)\Bigg) \\
&=: \Lambda_1 f + \Lambda_2 f
\end{aligned}
$$
where $K$ is a large positive parameter to be fixed later. Notice that in $\Lambda_1$, we have a term which is going to provide us some regularization
$$
\int_{\R^3 \times \mathbb{S}^2} B_\delta(v-v_*,\sigma) (\mu_* f' - \mu'_*f) \, \d\sigma \, \dv_*
$$
and another one which provides us some dissipativity:
$$
- K \langle v \rangle^{\gamma+2s} f.
$$

\subsection{Study of the controlled part $\Lambda_2$}
We first study the ``nice'' part of our splitting, namely $\Lambda_2$ which is ``almost bounded'' in the sense that it does not induce a loss of regularity but only a loss in weight.

\begin{lem} \label{lem:lambda_2}
Let $m(v):= \langle v \rangle^k$ with $k\ge0$. For any $K>0$ and for any $\ell>3/2$, we have the following estimate:
\begin{equation} \label{eq:lambda_2}
\|\Lambda_2 f\|_{H^\varsigma_{x,v}(m)} \lesssim \|f\|_{H^\varsigma_{x,v} (\langle v \rangle^{\gamma + 1 + \ell} m)}, \quad \forall \, \varsigma \in \R^+.
\end{equation}
\end{lem}
\preuve
We only look at the case $\varsigma \in \N$ and conclude that the result also holds for~$\varsigma \in \R^+$ by an interpolation argument. Let us begin with the case $\varsigma=0$ i.e. the $L^2$-case. We have
$$
\begin{aligned}
\Lambda_2 f &= K \langle v \rangle^{\gamma+2s} f + \int_{\R^3 \times \mathbb{S}^2} B_\delta(v-v_*,\sigma) (\mu'_* - \mu_*)f' \, \d\sigma \, \dv_* \\
&\quad +  \int_{\R^3 \times \mathbb{S}^2} B_\delta(v-v_*,\sigma) (\mu'_* - \mu_*)\, \d\sigma \, \dv_* \, f + \int_{\R^3\times \mathbb{S}^2} B_\delta^c(v-v_*,\sigma) \mu'_* f' \, \d\sigma \, \dv_*  \\
&\quad -  \int_{\R^3\times \mathbb{S}^2} B_\delta^c(v-v_*,\sigma) \mu_*\, \d\sigma \, \dv_* \, f + Q(f,\mu)\\
&=: \Lambda_{21}f + \Lambda_{22} f + \Lambda_{23}f + \Lambda_{24} f + \Lambda_{25} f + \Lambda_{26} f.
\end{aligned}
$$
The estimate on $\Lambda_{21}$ is obvious:
$$
\|\Lambda_{21}f\|_{L^2_{x,v}(m)} \lesssim \|f\|_{L^2_{x,v} (\langle v \rangle^{\gamma+2s}m)}.
$$
The analysis of $\Lambda_{23}$ is also easy to perform using the cancellation lemma from~\cite{ADVW-ARMA}, we have:
$$
\Lambda_{23}f = (S * \mu) f
$$
with $S$ satisfying the estimate $|S(z)| \lesssim|z|^{\gamma}$. We deduce that $|S * \mu|(v) \lesssim \langle v \rangle^{\gamma}$ and thus
$$
\|\Lambda_{23}f\|_{L^2_{x,v}(m)} \lesssim \|f\|_{L^2_{x,v} (\langle v \rangle^{\gamma}m)}.
$$
To treat $\Lambda_{24}$ and $\Lambda_{25}$, we use the fact that the kernel $B_\delta^c$ is not singular since the grazing collisions are removed. Since $|v'-v| \sim |v-v_*| \sin(\theta/2)$, we have:
$$
|B_\delta^c(v-v_*,\sigma)| \le b(\cos\theta) |v-v_*|^\gamma \mathds{1}_{|v'-v| \ge \delta} \lesssim  b(\cos\theta) |v-v_*|^{\gamma+1} \sin (\theta/2).
$$
Consequently, we obtain using that $m \lesssim m' m'_*$ that for $\ell>3/2$:
$$
\begin{aligned}
&\|\Lambda_{24}f\|^2_{L^2_{x,v}(m)} \lesssim \int_{\T^3 \times \R^3} \left(\int_{\R^3\times \mathbb{S}^2} b(\cos \theta) \sin (\theta/2)|v-v_*|^{\gamma+1} \mu'_* \, f' \,  \d\sigma \, \dv_* \right)^2 m^2 \, \dv \, \d x \\
&\lesssim \int_{\T^3 \times \R^3 \times \R^3 \times \mathbb{S}^2} b(\cos \theta) \sin (\theta/2) |v-v_*|^{2(\gamma+1)} (\mu'_* m'_*)^2 \, (f'm')^2 \, \langle v_* \rangle^{2\ell} \, \d\sigma \, \dv_* \,  \dv \, \d x
\end{aligned}
$$
where we have used Jensen inequality with the finite measure $b(\cos \theta) \sin\left(\theta / 2\right)d\sigma$ and Cauchy-Schwarz inequality with the measure $\langle v_* \rangle^\ell \dv_*$. Then, using the basic inequality
$\langle v_* \rangle \lesssim \langle v' \rangle \langle v'_* \rangle$ and the pre-post collisional change of variable, we get:
$$
\begin{aligned}
&\quad \|\Lambda_{24}f\|^2_{L^2_{x,v}(m)} \\
&\lesssim
\int_{\T^3 \times \R^3 \times \R^3 \times \mathbb{S}^2} b(\cos \theta) \sin (\theta/2) |v-v_*|^{2(\gamma+1)} (\mu_*m_*)^2 \, (fm)^2 \, \langle v \rangle^{2 \ell} \, \langle v_* \rangle^{2\ell} \, \d\sigma \, \dv_* \,  \dv \, \d x \\
&\lesssim
\|f\|^2_{L^2_{x,v} (\langle v \rangle^{\gamma+1+\ell} m)} \quad \text{with} \quad \ell>3/2.
\end{aligned}
$$
The treatment of $\Lambda_{25}$ is easier and we directly obtain:
$$
\|\Lambda_{25}f\|_{L^2_{x,v}(m)} \lesssim \|f\|_{L^2_{x,v} (\langle v \rangle^{\gamma+1} m)}.
$$
Concerning $\Lambda_{26}$, we have for any $\ell>3/2$:
$$
\|Q(f,\mu)\|_{L^2_v(m)} \lesssim  \|f\|_{L^2_v(\langle v \rangle^{\gamma+2s + \ell}m)}
$$
where we used~\cite[Theorem~2.1]{AMUXY-ARMA1}. We deduce that
$$
\|\Lambda_{26}f\|_{L^2_{x,v}(m)} \lesssim \|f\|_{L^2_{x,v} (\langle v \rangle^{\gamma+2s+\ell}m)}, \quad \ell >3/2.
$$
It now remains to deal with $\Lambda_{22}$. Denoting $M:= \sqrt{\mu}$, we have:
$$
\begin{aligned}
|\Lambda_{22} f| &\le \int_{\R^3 \times \mathbb{S}^2} B_\delta(v-v_*,\sigma) |M'_*-M_*| (M'_* + M_*) |f'| \, \d\sigma \, \dv_* \\
&\lesssim  \int_{\R^3 \times \mathbb{S}^2} b(\cos \theta) \sin (\theta/2) |v-v_*|^{\gamma+1} (M'_* + M_*) |f'| \, \d\sigma \, \dv_*
\end{aligned}
$$
where we used that the gradient of $M$ is bounded on $\R^d$. Then we use that $m \lesssim m' m'_*$ and $m \lesssim  \langle v-v_* \rangle^k m_*$ to get:
$$
\begin{aligned}
\|\Lambda_{22} f\|^2_{L^2_{x,v}(m)} &\lesssim \int_{\T^3 \times \R^3} \left(  \int_{\R^3 \times \mathbb{S}^2} b(\cos \theta) \sin (\theta/2) |v-v_*|^{\gamma+1} M'_*|f'| \, \d\sigma \, \dv_* \right)^2 m^2 \dv \, \d x\\
&\quad+  \int_{\T^3 \times \R^3} \left(  \int_{\R^3 \times \mathbb{S}^2} b(\cos \theta) \sin (\theta/2) |v-v_*|^{\gamma+1} M_*|f'| \, \d\sigma \, \dv_* \right)^2 m^2 \dv \, \d x\\
&\lesssim \int_{\T^3 \times \R^3} \left(  \int_{\R^3 \times \mathbb{S}^2} b(\cos \theta) \sin (\theta/2) |v-v_*|^{\gamma+1} M'_* m'_* |f'| m' \, \d\sigma \, \dv_* \right)^2 \dv \, \d x \\
&\quad + \int_{\T^3 \times \R^3} \left(  \int_{\R^3 \times \mathbb{S}^2} b(\cos \theta) \sin (\theta/2) \langle v-v_* \rangle^{\gamma+1+k} M_* m_*|f'| \, \d\sigma \, \dv_* \right)^2 \dv \, \d x \\
&=: I_1 + I_2.
\end{aligned}
$$
Using Jensen inequality and H\"{o}lder inequality as previously, we obtain for $\ell>3/2$:
$$
\begin{aligned}
&\|\Lambda_{22} f\|^2_{L^2_{x,v}(m)} \\
&\quad\lesssim \int_{\T^3 \times\R^3\times\R^3 \times \mathbb{S}^2} b(\cos \theta) \sin (\theta/2) |v-v_*|^{2(\gamma+1)} \mu'_* (m'_*)^2|f'|^2 (m')^2 \langle v_* \rangle^{2\ell} \, \d\sigma \, \dv_* \, \dv \, \d x\\
&\qquad+  \int_{\T^3 \times\R^3\times\R^3 \times \mathbb{S}^2} b(\cos \theta) \sin (\theta/2) \langle v-v_* \rangle^{2(\gamma+1+k)} \mu_* m_*^2 |f'|^2  \langle v_* \rangle^{2\ell} \,  \d\sigma \, \dv_* \, \dv \, \d x\\
&\quad =: I_1 + I_2.
\end{aligned}
$$
The first term $I_1$ is treated as $\Lambda_{24}$ and we thus have:
$$
I_1 \lesssim \|f\|^2_{L^2_{x,v} (\langle v \rangle^{\gamma+1+\ell} m)} .
$$
Concerning $I_2$, we first look at the integral
$$
J :=  \int_{\R^3 \times \mathbb{S}^2} b(\cos \theta) \sin (\theta/2)  \langle v-v_* \rangle^{2(\gamma+1+k)} |f'|^2 \, \d\sigma \, \dv.
$$
Then, for each $\sigma$, with $v_*$ still fixed, we perform the change of variables $v \rightarrow v'$. This change of variables is well-defined on the set $\left\{ \cos \theta > 0 \right\}$. Its Jacobian determinant is
$$
\left| \frac{dv'}{dv} \right| = \frac{1}{8} (1 + \kappa \cdot \sigma) = \frac{(\kappa' \cdot \sigma)^2}{4},
$$
where $\kappa=(v-v_*)/|v-v_*|$ and $\kappa'=(v'-v_*)/|v'-v_*|$. We have $\kappa' \cdot \sigma = \cos(\theta/2) \geq 1/\sqrt{2}$. The inverse transformation $v' \rightarrow \psi_\sigma(v')=v$ is then defined accordingly. We have
$$
\cos \theta = \kappa \cdot \sigma = 2 (\kappa' \cdot \sigma)^2 -1 \quad \text{and} \quad \sin (\theta/2) = \sqrt{ 1 - \cos^2(\theta/2)} = \sqrt{1 - (\kappa' \cdot \sigma)^2},
$$
and also
$$
|\psi_\sigma(v)-v_*|=|v-v_*|/\kappa\cdot\sigma.
$$
As a result, we get:
\begin{align*}
J&= \int_{\R^3 \times \mathbb{S}^2} b(2 (\kappa' \cdot \sigma)^2 -1) \sqrt{1 - (\kappa' \cdot \sigma)^2}  \langle \psi_\sigma(v')-v_* \rangle^{2(\gamma+1+k)} |f'|^2 \, \d\sigma \, \dv \\
&= \int _{ \kappa' \cdot \sigma  \geq 1/\sqrt{2}} b(2 (\kappa' \cdot \sigma)^2 -1) \sqrt{1 - (\kappa' \cdot \sigma)^2}  \langle \psi_\sigma(v')-v_* \rangle^{2(\gamma+1+k)} |f'|^2 \, \d\sigma \, \frac {4 \,\dv'}{(\kappa' \cdot \sigma)^2}\\
&= \int _{ \kappa \cdot \sigma  \geq 1/\sqrt{2}} b(2 (\kappa \cdot \sigma)^2 -1) \sqrt{1 - (\kappa \cdot \sigma)^2}  \langle \psi_\sigma(v)-v_* \rangle^{2(\gamma+1+k)} |f|^2 \, \d\sigma \, \frac {4 \,\dv}{(\kappa \cdot \sigma)^2} \\
&\lesssim \int _{ \kappa \cdot \sigma  \geq 1/\sqrt{2}} b(2 (\kappa \cdot \sigma)^2 -1) \sqrt{1 - (\kappa \cdot \sigma)^2}  \langle v-v_* \rangle^{2(\gamma+1+k)} |f|^2 \, \d\sigma \, \dv \\
&\lesssim \int_{\mathbb{S}^2} b(\cos (2\theta)) \sin(\theta) \, \d\sigma \int_{\R^3} f^2 m^2 \langle v \rangle^{2(\gamma+1)} \, \dv \, \langle v_* \rangle^{2(\gamma+1)} m_*^2.
\end{align*}
From this, we deduce that
$$
I_2 \lesssim  \|f\|^2_{L^2_{x,v} (\langle v \rangle^{\gamma+1} m)}
$$
and this concludes the proof in the case $\varsigma=0$.

Let us now explain briefly how to treat higher order derivatives: we only deal with the $H^1$-case, the other cases being handled similarly. For the derivative in $x$, we have immediately that for any $\ell>3/2$,
$$
\|\nabla_x \Lambda_2 f\|_{L^2_{x,v}(m)} \lesssim \|\nabla_xf\|_{L^2_{x,v} (\langle v \rangle^{\gamma + 1 + \ell} m)}
$$
since the operators $\nabla_x$ and $\Lambda_2$ commute ($\Lambda_2$ acts only in velocity). Concerning the derivative in $v$, we have to be more careful and in what follows, we only give the key points to obtain the final estimate. For the first term, we have:
$$
|\nabla_v \Lambda_{21} f| \lesssim \langle v \rangle^{\gamma+2s-1} |f| + \langle v \rangle^{\gamma+2s} |\nabla_v f|.
$$
For $\Lambda_{23}$, using the cancellation lemma, we have
$$
\nabla_v (\Lambda_{23} f) = (S*\nabla_v \mu) f + (S*\mu) \nabla_v f
$$
and we also have $|S*\nabla_v \mu| \lesssim \langle v \rangle^\gamma$.
For $\Lambda_{26}$ we can use the classical result (see~\cite{Vill-JMPA}) that tells us
$$
\nabla_v Q(f,\mu) = Q(\nabla_v f, \mu) + Q(f,\nabla_v \mu).
$$
In the same spirit that the latter formula is proven, one can show that
$$
\nabla_v \Lambda_{22}f =  \int_{\R^3 \times \mathbb{S}^2} B_\delta(v-v_*,\sigma) ((\nabla_v\mu)'_* - (\nabla_v\mu)_*)f' \, \d\sigma \, \dv_* + \Lambda_{22} (\nabla_vf),
$$
$$
\nabla_v \Lambda_{24} f = \int_{\R^3\times \mathbb{S}^2} B_\delta^c(v-v_*,\sigma) (\nabla_v\mu)'_* f' \, \d\sigma \, \dv_*   + \Lambda_{24} (\nabla_v f)
$$
and
$$
\nabla_v \Lambda_{25} f = -  \int_{\R^3\times \mathbb{S}^2} B_\delta^c(v-v_*,\sigma) (\nabla_v\mu)_*\, \d\sigma \, \dv_* \, f + \Lambda_{25}(\nabla_v f).
$$
The key elements to prove those relations are that $\nabla_v B_\delta = -\nabla_{v_*} B_\delta$ and that we have for any suitable function $g$:
$$
(\nabla_v + \nabla_{v_*})(g')=\left(\nabla_v g\right)' \quad \text{and} \quad (\nabla_v + \nabla_{v_*})(g'_*)=\left(\nabla_v g\right)'_*.
$$
Gathering the previous remarks, we are then able to obtain that for any $\ell>3/2$:
$$
\|\nabla_v \Lambda_2 f\|_{L^2_{x,v}(m)} \lesssim \|f\|_{L^2_{x,v} (\langle v \rangle^{\gamma + 1 + \ell} m)} + \|\nabla_vf\|_{L^2_{x,v} (\langle v \rangle^{\gamma + 1 + \ell} m)},
$$
which allows us to conclude.
\fin   \Bk

\subsection{Regularization properties of $\Lambda_1$ {in the dual case}} \label{subsec:reg}
The main result of this Subsection is Proposition~\ref{prop:lambda1*} and is about the regularization features of the semigroup associated to $\Lambda_1$. Here, we just state the result and we postpone its proof to Section~\ref{sec:quantization} in which we develop pseudo-differential arguments.

\subsubsection*{Functional spaces} In the remainder part of this section, we consider three weights:
$$
\left\{
\begin{aligned}
&\text{$m(v) = \langle v \rangle^k$ with {$k \ge 0$},}\\
&\text{$m_0(v) = \langle v \rangle^{k_0}$ with {$k_0 >\gamma/2 +3+2s$}}\\
&\text{$m_1(v) = \langle v \rangle^{k_1}$ with {$k_1=k_0 +\gamma + 1 + \ell$ and $\ell>3/2$}.}
\end{aligned}
\right.
$$
We then denote for $i=\emptyset,0,1$:
$$
\left\{
\begin{aligned}
&X_i=L^2_{x,v}(m_i)\\
&Y_i = H^{s,0}_{x,v}({ \langle v \rangle^{\gamma/2}} m_i)\\
&Z_i=H^{0,s}_{x,v}( \langle v \rangle^{\gamma/2} m_i)) \cap  L^2_{x,v}( \langle v \rangle^{(\gamma+2s)/2} m_i) \\
&\text{$Y'_i$ the dual of $Y_i$ w.r.t. $X_i$}\\
&\text{$Z'_i$ the dual of $Z_i$ w.r.t. $X_i$.}
\end{aligned}
\right.
$$
We also introduce the (almost) flat spaces:
$$
\left\{
\begin{aligned}
&F=L^2_{x,v} \\
&G=H^{s,0}_{x,v}({ \langle v \rangle^{\gamma/2}}) \\
&H=H^{0,s}_{x,v}( \langle v \rangle^{\gamma/2}) \cap  L^2_{x,v}( \langle v \rangle^{(\gamma+2s)/2}) \\
&\text{$G'$ the dual of $G$ w.r.t. $F$}\\
&\text{$H'$ the dual of $H$ w.r.t. $F$.}
\end{aligned}
\right.
$$

\subsubsection*{Remark on the dual embeddings}
First, we notice that
\begin{equation} \label{prop:embed}
\forall \, q_1 \le q_2, \, \varsigma \in \R^+, \quad H^\varsigma_v(\langle v \rangle^{q_2}) \hookrightarrow H^\varsigma_v(\langle v \rangle^{q_1}).
\end{equation}
This property is clear in the case $\varsigma \in \N$. Let us now treat the case $\varsigma \in \R^+ \setminus \N$. 
Since the weighted space $H^\varsigma_v(\langle v \rangle^{q_i})$ is defined through
$$
h \in H^\varsigma_v(\langle v \rangle^{q_i}) \Leftrightarrow h \langle v \rangle^{q_i} \in H^\varsigma_v
$$
and that we have, using the standard real interpolation notations (see for example~\cite{BookBL}):
$$
H^\varsigma_v = \left[H^{\lfloor \varsigma \rfloor}_v,H^{\lfloor \varsigma \rfloor+1}_v\right]_{{\varsigma-\lfloor \varsigma \rfloor},2},
$$
one can prove that
$$
H^\varsigma_v(\langle v \rangle^{q_i}) = \left[H^{\lfloor \varsigma \rfloor}_v(\langle v \rangle^{q_i}),H^{\lfloor \varsigma \rfloor+1}_v(\langle v \rangle^{q_i})\right]_{{\varsigma-\lfloor \varsigma \rfloor},2}, \quad i=1,2.
$$
From this, since $H^{\ell}_v(\langle v \rangle^{q_2}) \hookrightarrow H^{\ell}_v(\langle v \rangle^{q_1})$ for $\ell \in \N$, we deduce the desired embedding result:
$
H^\varsigma_v(\langle v \rangle^{q_2}) \hookrightarrow H^\varsigma_v(\langle v \rangle^{q_1}).
$

We can now prove that the standard inclusions for dual spaces do not hold here. Indeed, we have for example $Y_1 \subset Y_0$ and also $Y'_1 \subset Y'_0$ (the same for ``$Z$-spaces'' hold). This is due to the fact that the pivot spaces are $X_i$ and not $L^2_{x,v}$ as usually. Indeed, using that $k_1 \ge k_0$ and~\eqref{prop:embed}, we have
\begin{align*}
\| f\|_{Y'_0} &= \sup_{\|\varphi\|_{Y_0} \le 1} \langle f, \phi \rangle_{X_0} \\
&=  \sup_{\|\varphi m_0\|_{G} \le 1} \bigg\langle f m_1, \phi {m_0^2 \over m_1} \bigg\rangle_{F} \\
&= \sup_{\|\psi m_1^2/m_0\|_{G} \le 1} \langle f m_1, \psi m_1 \rangle_{F} \\
&\le  \sup_{\|\psi m_1\|_{G} \le 1} \langle f m_1, \psi m_1 \rangle_{F} = \sup_{\|\phi\|_{Y_1} \le 1} \langle f, \phi \rangle_{X_1} = \|f\|_{Y'_1}.
\end{align*}

\subsubsection*{Reduction of the problem to a ``simpler'' framework}
We start by explaining how to avoid some difficulties coming from the spaces in which we are working.
First, in order to simplify the problem, since we work in weighted spaces, we are going to ``include'' the weight in our operator.
For this purpose, we define the operator $\Lambda^{m}_1$ by
$$
\Lambda^{m}_1 f := m \, \Lambda_1 ({m}^{-1} f).
$$
We notice that if $f$ satisfies $\partial_t f = \Lambda_1 f$, then $h := {m}f$ satisfies $\partial_t h = \Lambda^{m}_1 h$ and we thus have $S_{\Lambda_1^{m}}(t) h = m S_{\Lambda_1}(t) f$.
Then, in order to avoid having to work in dual spaces, we introduce formal dual operators for which we prove regularization properties in ``positive'' Sobolev spaces.
To this end, we introduce the (formal) adjoint operator (w.r.t. the scalar product of $L^2_{x,v}$) of $\Lambda_1^m$ that we denote $\Lambda^{{m},*}_{1}$  and which is defined by:
\begin{equation*} 
\Lambda_{1}^{{m},*} \phi := \int_{\R^3 \times \mathbb{S}^2} B_\delta(v-v_*,\sigma) \, \mu'_* \, ( \phi' {m}'-\phi {m}) \, \d\sigma \, \dv_*\, {m}^{-1}  - K \langle v \rangle^{\gamma+2s} \, \phi + v \cdot \nabla_x \phi.
\end{equation*}
The advantage of working with this operator is that we can work in flat and positive Sobolev spaces.
We now write our main regularization estimate:

\begin{prop} \label{prop:lambda1*}
For $K$ large enough, we have the following estimates:
\begin{equation} \label{eq:dual}
\forall \, t \in (0,1], \quad \|S_{\Lambda^{{m},*}_1}(t)\phi\|_{H} \lesssim {1 \over \sqrt{t}} \|\phi\|_{F} \quad \text{and} \quad  \|S_{\Lambda^{{m},*}_1}(t)\phi\|_{G} \lesssim {1 \over t^{1/2+s}} \|\phi\|_{F}.
\end{equation}
\end{prop}
The proof of Proposition~\ref{prop:lambda1*} is to be compared with the one of Theorem~\ref{FKEthm}. Indeed, it is the same proof strategy, we introduce a functional which is going to be an entropy for our equation for small times. However, it is much more complicated in this case and our approach requires refined pseudo-differential tools, Section~\ref{sec:quantization} is dedicated to its proof. Before that, we explain how to use Proposition~\ref{prop:lambda1*} to get our final result in Theorem~\ref{thm:LIBthm}.

%


\subsection{Proof {of the dual result} of Theorem~\ref{thm:LIBthm}}
The goal is first to prove the dual result in Theorem~\ref{thm:LIBthm} in the case $r=0$. As already mentioned, in the case of the fractional Kolmogorov equation, the proof will be exactly the same for other values of $r$ since the operator $\Lambda^r_x$ commutes with the Boltzmann operator. We can thus apply the result obtained for $r=0$ to $\Lambda^r_x f_0$ to recover the result for $r\neq 0$.

From Proposition~\ref{prop:lambda1*}, we can deduce an estimate on the semigroup associated to $\Lambda_1$ in the ``original'' (non flat) spaces:

\begin{cor}
For $K$ large enough, the following estimates hold:
\begin{equation} \label{eq:lambda_1}
\forall \, t \in (0,1], \quad \|S_{\Lambda_1}(t)f\|_{X} \lesssim {1 \over \sqrt{t}} \|f\|_{Z'} \quad \text{and} \quad \|S_{\Lambda_1}(t)f\|_{X} \lesssim {1 \over t^{1/2+s}} \|f\|_{Y'}.
\end{equation}
\end{cor}

\preuve
Let us consider $K$ large enough so that the conclusion of Proposition~\ref{prop:lambda1*} holds. Using~\eqref{eq:dual}, we have for any $t \in (0,1]$:
$$
\begin{aligned}
\|S_{\Lambda_1}(t) f \|_{X} &\lesssim \|S_{\Lambda_1^{m}} (t) h \|_F = \sup_{\|\phi\|_{F} \le 1} \langle S_{\Lambda_1^m}(t) h, \phi \rangle = \sup_{\|\phi\|_{F} \le 1} \langle h, S_{\Lambda_1^{{m},*}}(t)\phi \rangle \\
&\lesssim \sup_{\|\phi\|_{F} \le 1} \|h\|_{H'} \|S_{\Lambda_1^{{m},*}}(t)\phi\|_{H} \lesssim  {1 \over \sqrt{t}}\|h\|_{H'} \lesssim {1 \over \sqrt{t}} \|f\|_{Z'}
\end{aligned}
$$
which is exactly the first part of~\eqref{eq:lambda_1}. The second one is proven in the same way.
\fin

Let us finally prove that the regularization properties of $\Lambda_1$ are enough to conclude that the whole operator $\Lambda$ has some good regularization properties: even if we have a loss of weight in the final estimate, $\Lambda$ inherits regularization properties from $\Lambda_1$ in terms of fractional Sobolev norms.

\begin{lem} \label{lem:lambda_1/lambda}
We have:
\begin{equation} \label{eq:lambda}
\forall \, t \in (0,1], \quad \|S_{\Lambda}(t)f\|_{X_0} \lesssim {1 \over \sqrt{t}} \|f\|_{Z'_1} \quad \text{and} \quad \|S_{\Lambda}(t)f\|_{X_0} \lesssim {1 \over t^{1/2+s}} \|f\|_{Y'_1}.
\end{equation}
\end{lem}

\preuve
We recall that from~\cite{HTT1*}, we have that $\Lambda$ generates a semigroup in $X_0$ and thus we have the estimate
\begin{equation} \label{eq:lambdaHTT1*}
\forall \, t \in (0,1], \quad \|S_{\Lambda}(t) f \|_{X_0} \lesssim \|f\|_{X_0}.
\end{equation}
Then, we write Duhamel formula:
$$
S_{\Lambda}(t) = S_{\Lambda_1}(t) + \int_0^t S_{\Lambda}(s) \Lambda_2 S_{\Lambda_1}(t-s) \, \d s
$$
from which we deduce, combining~\eqref{eq:lambdaHTT1*},~\eqref{eq:lambda_1} and~\eqref{eq:lambda_2} applied with the appropriate weights, that for $t \in (0,1]$,
$$
\begin{aligned}
\|S_{\Lambda}(t)h\|_{X_0}
&\lesssim \|S_{\Lambda_1}(t) f\|_{X_0} + \int_0^t \|S_{\Lambda}(s)\Lambda_2 S_{\Lambda_1}(t-s) f\|_{X_0} \, \d s \\
&\lesssim {1 \over \sqrt{t}} \|f\|_{Z'_0} + \int_0^t \|\Lambda_2 S_{\Lambda_1}(t-s) f\|_{X_0} \, \d s \\
&\lesssim {1 \over \sqrt{t}} \|f\|_{Z'_0} + \int_0^t \|S_{\Lambda_1}(t-s) f\|_{X_1} \, \d s \\
&\lesssim {1 \over \sqrt{t}} \|f\|_{Z'_0} + \int_0^t {1 \over \sqrt{t-s}} \|f\|_{Z'_1} \, \d s \\
&\lesssim {1 \over \sqrt{t}} \|f\|_{Z'_0} + \int_0^1 {1 \over \sqrt{s}} \|f\|_{Z'_1} \, \d s \lesssim {1 \over \sqrt{t}} \|f\|_{Z'_1}.
\end{aligned}
$$
This concludes the proof of the first part of~\eqref{eq:lambda}. Concerning the second one, we proceed as before using that $1/2+s<1$ since $s<1/2$ and we obtain for any $t \in (0,1]$:
$$
\begin{aligned}
\|S_{\Lambda}(t)h\|_{X_0}
\lesssim {1 \over t^{1/2+s}} \|f\|_{Y'_1}.
\end{aligned}
$$
\fin

\section{Proof of Proposition \ref{prop:lambda1*}} \label{sec:quantization}

The aim of this section is the proof of Proposition \ref{prop:lambda1*} about the regularization properties of the operator
\begin{equation*} 
\Lambda_{1}^{{m},*} \phi = \int_{\R^3 \times \mathbb{S}^2} B_\delta(v-v_*,\sigma) \, \mu'_* \, ( \phi' {m}'-\phi {m}) \, \d\sigma \, \dv_*\, {m}^{-1}  - K \langle v \rangle^{\gamma+2s} \, \phi + v \cdot \nabla_x \phi.
\end{equation*}
This will be done with a pseudodifferential version of the Lyapunov trick developed in the Fokker-Planck case and special classes of symbols that we recall in the Appendix.

\subsection{Pseudodifferential formulation of the operator $\Lambda_1^{m,*}$}

The operator $\Lambda_{1}^{{m},*}$ is very similar to the operator $\mathcal{L}_{1,2,\delta}$ defined in \cite[Proposition 3.1]{AHL*}. We shall thus take advantage of the analysis of the pseudo-differential operator $\mathcal{L}_{1,2,\delta}$ and its symbol in \cite{AHL*}.
If we extract the collision part of the operator $\Lambda_1^{m,*}$ (forgetting the transport one and the addition of the multiplicative term), we obtain
$$
\Lambda_{1}^{{m},*, {\hbox{\footnotesize collision}}} \phi := \int_{\R^3 \times \mathbb{S}^2} B_\delta(v-v_*,\sigma) \, \mu'_* \, ( \phi' {m}'-\phi {m}) \, \d\sigma \, \dv_*\, {m}^{-1}
$$
In the case $m= 1$, this operator is actually the main one studied in~\cite{AHL*}:
$$
\Lambda_{1}^{{1},*, {\hbox{\footnotesize collision}}} = \mathcal{L}_{1,2,\delta} =: -\tilde{a}_0(v,D_v) ,
$$
where $\tilde{a}_0$ is a real symbol in $(v,\eta)$ defined through
$$
\tilde{a}_0(v,\eta) := \int_{\R^3_h} {{\d h} \over |h|^{3+2s}}  \int_{E_{0,h}} \d\alpha \, \widetilde{b}(\alpha,h) \, \mathds{1}_{|\alpha| \ge |h|} \, \chi_\delta(h) \, \mu(\alpha+v) \, |\alpha+h|^{\gamma+1+2s} \,
(1-\cos(\eta \cdot h))
$$
thanks to Carleman representation (see Lemma~\ref{lem:Carleman}).  We  recall here the main result from~\cite{AHL*} concerning the symbol $\tilde{a}_0$ (be careful, this symbol is denoted without tilde there):

\begin{prop}[Propositions~3.1 and~3.4 in \cite{AHL*}] \label{prop:AHL}
The symbol $\tilde{a}_0$ satisfies the following properties:
\begin{equation*}
\begin{split}
i) &  \qquad  \tilde{a}_0 \in S( \seq{v}^\gamma (1+ |\eta|^2 + |v\wedge\eta|^2)^s, \Gamma), \\
 ii) &  \qquad  \forall \,\eps>0, \, \,  \nabla_\eta \tilde{a}_0 \in S( \eps \seq{v}^\gamma (1+ |\eta|^2 + |v\wedge\eta|^2)^s + \eps^{-1} \seq{v}^{\gamma + 2s}), \\
iii)  & \qquad   \exists \, c>0, \, -c\seq{v}^{\gamma + 2s}+ \seq{v}^\gamma \sep{ 1+ |\eta|^2 + |v\wedge\eta|^2}^s \lesssim \tilde{a}_0 \lesssim \seq{v}^\gamma \sep{ 1+ |\eta|^2 + |v\wedge\eta|^2}^s,
\end{split}
\end{equation*}
where $\Gamma:=|\dv|^2+|\d \eta|^2$ is the flat metric.
\end{prop}

For convenience we denote by $a_0$ the Weyl symbol of operator $\tilde{a}_0(v,D_v)$, so that
$$
a_0^w = \tilde{a}_0(v,D_v).
$$
Everywhere in what follows, any symbol with a tilde will refer to a classical quantization, and when no tilde is present, the symbol will refer to the Weyl quantization. Both quantizations are recalled in the beginning of Section \ref{pseud} in the Appendix.
 Note that $a_0$ is not real anymore, anyway we shall see later that it conserves good ellipticity properties. Denoting then
$$
a(v,\eta) := \sep{  m^{-1} \sharp {a}_0 \sharp m } (v,\eta)+ K \seq{v}^{\gamma + 2s},
$$
where $\sharp$ denotes the usual Weyl composition and we omit the dependency of $a$ with respect to $K$ in our notation, we have:
$$
\Lambda_{1}^{{m},*} = -a^w + v \cdot \nabla_x.
$$
For sake of simplicity, we introduce the following notation
$$
A := a^w,
$$ so that  the collision part of
operator $\Lambda_{1}^{{m},*}$ writes
$$
 \Lambda_{1}^{{m},*} = -A + v \cdot \nabla_x
$$
(recall that they depend on $K$).
In order to study the symbolic properties of $a$, 
we now introduce
the main weights. We pose for $(v,\eta) \in \R^6$
$$
\lambda_v^2(v,\eta) :=  \seq{\eta}^2 + \seq{v\wedge \eta}^2 +  \seq{v}^2
$$
and for given $s \in (0,1/2)$ and $\gamma \in (0,1)$ we pose
$$
p(v,\eta) := \seq{v}^\gamma \lambda_v^{2s} + K \seq{v}^{\gamma + 2s}
$$
which will be the main reference symbol of our study (note that this symbol is denoted $\tilde{a}_K$ in~\cite{AHL*}). Although $p$ depends on $K$, we will omit in the following any subscript or reference to this dependence. It will be shown in the next subsection that $p$ is a good weight in the sense of the Appendix. The following Lemma shows that
$a$ has good properties in the class $S_K(p)$, the main class of symbols whose definition is recalled in full generality in the Appendix.

 \begin{lem} \label{aposell} Let $m(v) = \seq{v}^k$ for $k\in \R$. Then
 uniformly in $K$ sufficiently large, we have that  $\Re a \geq 0$,  $a \in S_K(p)$ and $\Re a$ is elliptic positive in this class. 
  \end{lem}
\preuve We shall take profit of the estimates from~\cite{AHL*} recalled above in Propostion~\ref{prop:AHL}. We first note that
because of the symbolic estimates on $\tilde{a}_0$ we can take $\eps= K^{-1/2}$ in ii) and,  using Lemma \ref{allsk}, we get that
$\tilde{a}_0 \in S_K(p)$ and then ${a}_0 \in S_K(p)$.
Adding $K \seq{v}^{\gamma + 2s}$ does not change the computation and we also get
that
$$
a_0 + K \seq{v}^{\gamma + 2s} \in S_K(p).
$$
Now we can do the conjugation with $m$. We first note that clearly, with the same notations as before, we have
$ m \in S_K(m)$ and $m^{-1} \in S_K(m^{-1})$.
This can be checked directly by noticing that the derivatives of $m$ in $\eta$ are zero. The stability of the class $S_K$ from Lemma~\ref{stab} implies then  that
$$
a  = m^{-1} \sharp  a_0   \sharp  m + K \seq{v}^{\gamma + 2s} = m^{-1} \sharp \sep{ a_0 + K \seq{v}^{\gamma + 2s} } \sharp  m \in S_K( p ).
$$
We can also notice that looking at the main terms in the asymptotic development of the~$\sharp$ product { (see in particular Lemma \ref{allsk} and its proof)},  we have
$$
a =  {a}_0  + K \seq{v}^{\gamma + 2s} + r =  \tilde{a}_0  + K \seq{v}^{\gamma + 2s} + r'
$$
 with  $r$ and $r' \in K^{-1/2} S(p)$ (note that $r$ is exactly  the { Weyl} symbol of $m^{-1} [ a_0^w ,  m ]$).
Since from Propostion~\ref{prop:AHL}.~iii), we have   $\tilde{a}_0 + K \seq{v}^{\gamma + 2s} \gtrsim p $ (uniformly in $K$), we get that
 $$
 \Re a \gtrsim p
 $$
 so that $\Re a$ is non-negative and elliptic for $K$ large  (note that this proof is very close to the one of Lemma \ref{allsk} in the appendix). \Bk
\fin

\subsection{Reference weights} \label{subsec:refweights}

We now introduce some weights involving the constant $K$ where~$K$ is a large constant to be defined later. Formally, $1/\sqrt{K}$ plays the role of a small semiclassical parameter.
We recall that for $(v,\eta) \in \R^6$
$$
\lambda_v^2(v,\eta) =  \seq{\eta}^2 + \seq{v\wedge \eta}^2 +  \seq{v}^2
$$
and for given $s \in (0,1/2)$ and $\gamma \in (0,1)$.
$$
p(v,\eta) = \seq{v}^\gamma \lambda_v^{2s} + K \seq{v}^{\gamma + 2s}.
$$
We shall need their counterparts in the $\xi$ variable (considered as a parameter) instead of~$\eta$ and thus also introduce
$$
\lambda_x^2(v,\eta) :=  \seq{\xi}^2 + \seq{v\wedge \xi}^2 +  \seq{v}^2
$$
and
$$
q(v,\eta) := \seq{v}^\gamma \lambda_x^{2s} + K \seq{v}^{\gamma + 2s},
$$
where we omit the dependance on $K$ and $\xi$ again in the notations. We eventually introduce a mixed symbol
$$
\omega (v,\eta): = { - }\langle v\rangle ^\gamma \lambda_x^{s-1} \lambda_v^{s-1}(\eta\cdot \xi+ (v\wedge \eta) \cdot (v\wedge \xi))
$$
which will be crucial in the analysis.

Following the Appendix, we have in particular:
\begin{lem}
The symbols $p$, $q$ and more generally $\seq{v}^\zeta p^\varrho q^\varsigma$ for $\zeta$, $\varrho$ and $\varsigma \in \R$ are temperate with respect to $\Gamma$ uniformly w.r.t. $K$ and $\xi$.
\end{lem}
\preuve These computations are done for e.g. in~\cite[Section 3.3]{AHL*}.
\fin

The symbols $p$, $q$, and $\omega$ are then good symbols w.r.t. these classes, as the  following lemma shows.

\begin{lem} \label{stabsymb} We have $p \in S_K(p)$, $q\in S_K(q)$, $\omega \in S_K(\sqrt{pq})$ and more generally $\seq{v}^\zeta p^\varrho q^\varsigma \in S_K(\seq{v}^\zeta p^\varrho q^\varsigma )$ for $\zeta$, $\varrho$ and $\varsigma \in \R$, all this uniformly in $K$ and $\xi$. \end{lem}

\preuve
We only do the proof for $p$, the other being similar. We just have to differentiate the symbol $p$. We study first the gradient with respect to~$\eta$ (which corresponds to the case $|\beta|=1$). We notice that
\begin{equation*}
\begin{split}
\nabla_\eta p & = s \seq{v}^\gamma \lambda_v^{2s-2} \nabla_\eta (\lambda_v^2).
\end{split}
\end{equation*}
We also have that
$$
\abs{\nabla_\eta (\lambda_v^2)}
 \leq 2 \lambda_v \seq{v}
$$
from which we deduce that
\begin{align*}
\abs{ \nabla_\eta p } &\le 2s \seq{v}^{\gamma+1} \lambda_v^{2s-1} \\
&= 2s K^{-1/2} \left(K^{1/2} \seq{v}^{\gamma/2+s}\right) \seq{v}^{\gamma/2+1-s} \lambda_v^{2s-1} \\
&\le 2s K^{-1/2} p^{1/2} \seq{v}^{\gamma/2} \lambda_v^{s} \le 2s K^{-1/2} p
\end{align*}
which is the desired result.
We skip the other similar computations. \Bk
\fin

\subsection{Technical lemmas} \label{subsec:technical}

The main idea in the proof of the regularization result in Proposition~\ref{prop:lambda1*} is to
mimic the proof of the Fokker-Planck case, using deeply the positivity preserving property of the Wick quantization.

In what follows, we state a series of lemmas (from~\ref{lemma1} to~\ref{lemma4}) which are crucial to be able to ``compare'' our operator $A$ with quantizations of the simpler symbols $p$ and $q$ we introduced in the preceding subsection.

\begin{lem}\label{lemma1}
There exists $c_a>0$ such that
$$
2 \Re \sep{ Af, f}\geq c_a
 \sep{ p ^{\rm Wick}f, f}.$$
\end{lem}

\preuve  We first notice that
$$
 \Re \sep{ Af, f} = \Re \sep{ a^w \Bk f,f}
= \sep{ ( \Re a)^w \Bk f,f}
$$
thanks to the properties of the Weyl quantization.
%
%
%
Using \eqref{mainscalar} for $\Re a$, we therefore get that
$$
 \Re \sep{ Af, f} = \sep{( \Re a)^w f, f} \simeq    \sep{ (\Re a)^\wick f, f} =
\Re  \sep{ a^\wick f, f}.
 $$
%
%
 Moreover, $\Re a \simeq p$ uniformly in $K$ from Lemma~\ref{aposell}. This implies that there exists $c_a >0$ such that $\Re a - c_a p \geq 0$. Using the positivity property of the Wick quantization gives  $\Re (a)^\wick \Bk - c_a p^\wick \geq 0$ in the sense of operators.
 %
%
  This proves the result.
 \fin
%
%
%
%

\begin{lem}\label{lemma2}
There exists $c_p>0$ such that
$$
\sep{p^{\rm Wick} A f + A^* p^{\rm Wick} f,f}   \geq c_p \ ((p^2)^{\rm Wick}f,f).
 $$
\end{lem}

\preuve 
We have from the definition of the Wick quantization (see \eqref{wickN2})
 $$
 p^\wick A + A^* p^\wick = \sep{ (p\star  N ) \sharp  a + \bar{a}  \sharp (p\star  N ) }^w.
 $$
Using now Lemma \ref{allsk}, we have that  $p \in S_K(p)$ implies $p\star  N  \in S_K(p)$ and
 from the second point in Lemma \ref{stab}, we get that
 $(p\star  N ) \sharp a + \bar{a} \sharp (p\star  N )$ is elliptic, real and positive (from selfadjointness) in $S_K(p^2)$. We therefore get from \eqref{mainscalar} that
 $$
 (\sep{ (p\star N) \sharp a + \bar{a} \sharp (p\star N) }^w f,f) \simeq
 (\sep{ (p\star N) \sharp a + \bar{a} \sharp (p\star N) }^\wick f,f)
$$
Since $(p\star N) \sharp a + \bar{a} \sharp (p\star N) \simeq p^2$ (uniformly in $K$), the positivity properties of the Wick quantization imply the result.
\fin

\begin{lem}\label{lemma3} There exists $c_q>0$  such that
$$
\sep{q^{\rm Wick} A f + A^* q^{\rm Wick} f,f}\geq c_q \sep{(p q)^{\rm Wick} f,f}.
$$
\end{lem}

\preuve The proof is almost the same as the one of Lemma \ref{lemma2}, the main difference being that
the symbol $q$ now depends on a parameter $\xi$,  with respect to which all estimates have to be uniform. We write
 $$
 q^\wick A + A^* q^\wick = \sep{ (q\star N) \sharp {a} + \bar{{a}} \sharp (q\star N) }^w
 $$
 where again ${a}$ denotes the Weyl symbol of $A$.
 We have that $q \in S_K(q)$ uniformly in $K$ and $\xi$ and this implies  $q\star N \in S_K(q)$. From  ${a} \in S_K(p)$
 and the second point in Proposition~\ref{stab}, we get that
 $(q\star N) \sharp {a} + \bar{{a}} \sharp (q\star N)$ is elliptic, real and positive  in $S_K(pq)$. Together with~\eqref{mainscalar},  this implies that there exists $c_q>0$ s.t.
 $$
 (\sep{ (q\star N) \sharp {a} + \bar{{a}} \sharp (q\star N) }^w f,f) \simeq
 (\sep{ (q\star N) \sharp {a} + \bar{{a}} \sharp (q\star N) }^\wick f,f) \geq c_q ( (pq)^\wick f, f)
$$
 where the last inequality comes from the positivity properties of the Wick quantization.
 \fin

\begin{lem} \label{lemma3bis} There exist  $c_\omega>0$ such that
$$
\abs{ \sep{\omega^{\rm Wick} A f + A^* \omega^{\rm Wick} f,f}} \leq c_\omega  \sep{ (p^{3/2} q^{1/2})^\wick f,f}.
$$
\end{lem}

\preuve
We begin by denoting $\theta := p^{3/4} q^{1/4}$ so that $\theta^2 = p^{3/2} q^{1/2}$. Using Lemma \ref{stabsymb}, we get that $\theta$ is elliptic positive in $S_K(\theta)$. Note also
 that
 $$
 \omega^{\rm Wick} A  + A^* \omega^{\rm Wick} = \sep{ (\omega\star N ) \sharp {a} + \bar{{a}} \sharp (\omega\star N ) }^w
 $$
 using the definitions of the Wick quantization and still denoting again ${a}$ the Weyl symbol of operator $A$. From Lemma \ref{stabsymb}, $\omega \in S_K(\sqrt{pq})$ so that $\omega\star N$ is also in $S_K(\sqrt{pq})$ by Lemma~\ref{allsk}. On the other hand, ${a} \in S_K(p)$ and using the stability Proposition \ref{stab}, we therefore get that
 \begin{equation} \label{theta2}
  (\omega\star N ) \sharp {a} + \bar{{a}} \sharp (\omega\star N ) \in S_K(p^{3/2} q^{1/2}) = S_K(\theta^2).
 \end{equation}
 We then write
 \begin{equation*}
 \begin{split}
 & \abs{ \sep{\omega^{\rm Wick} A f + A^* \omega^{\rm Wick} f,f}} \\
 & = \bigg| \bigg(
    \underbrace{  (\theta^{-1})^\wick \sep{ (\omega\star N ) \sharp {a} + \bar{{a}} \sharp (\omega\star N )}^w
    (\theta^{-1})^\wick}_{ \textrm{ Operator } \Omega } ((\theta^{-1})^\wick)^{-1} f, ((\theta^{-1})^\wick)^{-1} f \bigg) \bigg|.
    \end{split}
    \end{equation*}
Let us prove that operator $\Omega$ is bounded. For this, we first note that
  $(\theta^{-1})^\wick = (\theta^{-1} \star N)^w$  and recall that $\theta$ is elliptic positive. Lemma \ref{stab} implies that $\theta^{-1}$ is positive elliptic in $S_K(\theta^{-1})$ too and from Lemma \ref{allsk}, the same is true for $\theta^{-1} \star N$. The Weyl symbol of $\Omega$ can be written
  $$
  \hbox{symb}(\Omega) = (\theta^{-1} \star N) {\sharp} \sep{ (\omega\star N ) \sharp {a} + \bar{{a}} \sharp (\omega\star N )} {\sharp}(\theta^{-1} \star N)
  $$
    and from the stability Lemma \ref{stab} and \eqref{theta2}, this symbol is in  $S_K(1)$. In particular, the operator $\Omega$ is bounded on $L^2$. We have that
    \begin{equation} \label{secondcomplique}
    \begin{split}
    |(\Omega ((\theta^{-1})^\wick)^{-1}f,((\theta^{-1})^\wick)^{-1}f)| & \leq C \norm{  ((\theta^{-1})^\wick)^{-1} f}^2 \\
    &\leq C' \norm{  \theta^\wick f}^2 \leq C'' ( (\theta^2)^\wick f,f).
 \end{split}
 \end{equation}
 The first inequality comes from the fact that $\Omega$ is bounded. The last inequality is just a consequence of \eqref{mainnorm}. Let us precise the arguments used for proving the second inequality:
  we have
   \begin{equation} \label{trucdebut}
    \begin{split}
     \norm{  ((\theta^{-1})^\wick)^{-1} f}^2  =  \norm{  ((\theta^{-1}\star N)^w)^{-1} f}^2    \simeq \norm{  ((\theta^{-1}\star N)^{-1})^w f}^2
      \end{split}
      \end{equation}
      using the definition of the Wick quantization and  \eqref{maingamma}. We also check by direct computation that $(\theta^{-1}\star N)^{-1}$ is elliptic positive in
      in $S_K(\theta)$ using Lemmas \ref{allsk} (see also Remark~\ref{allskrem}) and \ref{stab} b). This implies by \eqref{mainequiv} applied with $\tau = (\theta^{-1}\star N)^{-1}$ that
      \begin{equation} \label{trucmilieu}
      \norm{  ((\theta^{-1}\star N)^{-1})^w f}^2  \simeq  \norm{  \theta^w f}^2,
      \end{equation}
      and we get then by \eqref{mainnorm}
      \begin{equation} \label{trucfin}
      \norm{  \theta^w f}^2 \simeq ( (\theta^2)^\wick f,f).
      \end{equation}
      The estimate (\ref{trucdebut}-\ref{trucfin}) yield the second inequality in \eqref{secondcomplique}.

\Bk

\fin

To conclude this subsection, we state a lemma which will be useful in the sequel, and whose proof is direct using positivity properties of the Wick quantization.
\begin{lem} \label{lemma4}
We have the following estimates:
 $$
  \sep{(\langle v\rangle^{2\gamma}\lambda_v^{4s})^{\rm Wick} f,f} \le  \sep{(p^2)^{\rm Wick} f,f} \leq 2(1+K^2)   \sep{(\langle v\rangle^{2\gamma}\lambda_v^{4s})^{\rm Wick} f,f} ,
 $$
 $$
 \sep{p^{\rm Wick} f,f}= \sep{(\langle v\rangle^{\gamma}\lambda_v^{2s})^{\rm Wick} f,f}+K\sep{(\langle v\rangle^{\gamma+2s})^{\rm Wick} f,f},
 $$
  $$
 \sep{(\langle v\rangle^{2\gamma}\lambda_v^{2s}\lambda_x^{2s})^{\rm Wick} f,f} \le  \sep{(pq)^{\rm Wick} f,f} \le (1+K)^2 \sep{(\langle v\rangle^{2\gamma}\lambda_v^{2s}\lambda_x^{2s})^{\rm Wick} f,f} .
 $$
\end{lem}

\subsection{The Lyapunov functional} 
\label{subsec:proplambda1*}
From now on, we fix once and for all the constant~$K$ so that the conclusions of Lemmas \ref{lemma1} to \ref{lemma4} are true.
In the same spirit as  in Subsection~\ref{sublyapKFP} for the Fokker-Planck case, we build below a Lyapunov functional corresponding to the following equation
$$
\partial_t \varphi = v \cdot \nabla_x \varphi - A \varphi,
$$
and we consider  $\varphi$ a solution.
Then, since $A$ acts only on the velocity variable, we can take the Fourier transform of our equation in $x \in \T^3$ and see the associated Fourier variable $\xi\in \Z^3$ as a parameter in our equation. We thus consider $\psi=\mathcal{F}_x \varphi$ to be a solution of
$$
\partial_t \psi { -} i v\cdot \xi \psi + A \psi =0
$$
 with initial data $\psi_0$. We then follow the lines of the proof given in Section \ref{sec:FK} and we introduce an adapted entropy
functional defined  for all $t\geq 0$ by
\begin{equation} \label{hhh1}
\hhh(t) := C\norm{\psi}^2 + D t \sep{p ^{\rm Wick} \psi,\psi}+ E t^{1+s}  \sep{ \omega^{\rm Wick} \psi,  \psi} + t^{1+2s} \sep{q^{\rm Wick} \psi,\psi}
\end{equation}
for large constants $C$, $D$, $E$ to be chosen later, where $\norm{\cdot}$ is the usual $L^2$ norm and $\sep{ \cdot,\cdot}$ is the usual (complex) $L^2$ scalar product.

\begin{lem} \label{hh1bis}
If $E \leq \sqrt{D}$ then for all  $t\geq 0$, we have $\hhh(t) \geq 0$. Precisely, we have
\begin{equation*}
0 \leq  C\norm{\psi}^2 + \frac{D}{2} t \sep{p ^{\rm Wick} \psi,\psi}  + \frac{1}{2} t^{1+2s} \sep{q^{\rm Wick} \psi,\psi}
\leq \hhh(t).
\end{equation*}
\end{lem}
\preuve The first part of the inequality comes from the positivity property~\eqref{wick2}. For the bound on $\mathcal{H}(t)$, we start by noticing that using Cauchy-Schwarz inequality:
$$
|\eta\cdot \xi+ (v\wedge \eta) \cdot (v\wedge \xi)| \le \lambda_x \lambda_v.
$$
Then,  the time-dependent Cauchy-Schwarz inequality gives
$$
- Et^{s} \langle v\rangle ^\gamma \lambda_x^{s-1} \lambda_v^{s-1}(\eta\cdot \xi+ (v\wedge \eta) \cdot (v\wedge \xi))
\le \frac{E^2}{2} \langle v \rangle^\gamma \lambda_v^{2s} + \frac{1}{2} t^{2s} \langle v \rangle^{\gamma} \lambda_x^{2s}.
$$
The positivity of the Wick quantization and the fact that $E^2 \le D$  imply that
$$
E t^{1+s}  \sep{ \omega^{\rm Wick} \psi,  \psi} \ge -{D \over 2} t \sep{p ^{\rm Wick} \psi,\psi} - {1 \over 2} t^{1+2s} \sep{q^{\rm Wick} \psi,\psi}
$$
which proves the statement.
\fin

We now show that $\hhh$ is indeed a Lyapunov function (entropy functional).

\begin{lem} \label{derivB}  For well chosen (arbitrarily large) constants $C$, $D$ and $E$, we have
$$
\ddt \hhh(t) \leq 0, \quad \forall \, t \in (0,1].
$$
\end{lem}

\preuve
Let us define
$$
\mathcal P := p ^{\rm Wick} A + A^* p ^{\rm Wick},\quad \mathcal {\it\Omega} := \omega^{\rm Wick} A + A^* \omega^{\rm Wick}, \quad  \mathcal Q := q^{\rm Wick} A + A^* q^{\rm Wick}.
$$
Then, for each term in the sum, we have
\begin{equation}\label{C}
\ddt C \norm{\psi}^2= - 2C \, \Re \sep{ A\psi, \psi},
\end{equation}

\begin{equation}\label{D}
\ddt \left(D t  \sep{p ^{\rm Wick} \psi,\psi}\right)= D\sep{p^{\rm Wick} \psi,\psi} -   Dt \sep{\mathcal P \psi,\psi}{ +}D t \sep {\{p,v\cdot \xi\}^{\rm Wick} \psi,\psi},
\end{equation}

\begin{equation}\label{E}
\begin{aligned}
&\ddt\left( E t^{1+s} \sep{ \omega^{\rm Wick} \psi,  \psi} \right) \\
&\qquad =  (1+s)Et^{s}\sep{\omega^{\rm Wick} \psi,\psi} - Et ^{1+s}\sep{{\it \Omega} \psi,\psi}{ +} Et^{1+s}\sep{\{\omega,v\cdot \xi\}^{\rm Wick} \psi,\psi},
\end{aligned}
\end{equation}

\begin{equation}\label{F}
\begin{aligned}
&\ddt \left( t^{1+2s} \sep{ q^{\rm Wick} \psi,  \psi} \right)\\
&\qquad =(1+2s) t^{2s}\sep{q^{\rm Wick} \psi,\psi} - t^{1+2s} \sep{\mathcal Q \psi,\psi}{ +} t^{1+2s} \sep {\{q,v\cdot \xi\}^{\rm Wick} \psi,\psi},
\end{aligned}
\end{equation}
where, in the first term we used the skew-adjointness of the transport operator and in the last term of \eqref{D}, \eqref{E}, \eqref{F}, we used \eqref{wick}.

The right hand side in \eqref{C} is non-positive (thanks to the property of positivity of the Wick quantization~\eqref{wick2}) and using Lemma \ref{lemma1} and Lemma \ref{lemma4}, it can be estimated as
\begin{align*}
-2 C\Re \sep{ A\psi, \psi}&\leq - c_a C
 \sep{ p^{\rm Wick} \psi, \psi}\\ &\leq -\underbrace{c_a C\sep{(\langle v\rangle^{\gamma}\lambda_v^{2s})^{\rm Wick} \psi,\psi}}_{I}-\underbrace{ c_a CK\sep{(\langle v\rangle^{\gamma+2s})^{\rm Wick} \psi,\psi}}_{II}.
 \end{align*}

 Analogously, we can  deduce a bound for the first term in \eqref{D}. Indeed, we recover two non-negative terms
 $$
 D\sep{ p^{\rm Wick} \psi, \psi}\leq \underbrace{D\sep{(\langle v\rangle^{\gamma}\lambda_v^{2s})^{\rm Wick} \psi,\psi}}_{i}+\underbrace{D K\sep{(\langle v\rangle^{\gamma+2s})^{\rm Wick} \psi,\psi}}_{ii}.
 $$
 Moreover,  using the positivity of the Wick quantization~\eqref{wick2}, the second  term in \eqref{D} is non-positive and, using Lemma \ref{lemma2} and  Lemma \ref{lemma4}, it can be estimated as
 $$
 -Dt\sep{\mathcal P  \psi,\psi}\leq -c_p D t\sep{(p ^2)^{\rm Wick} \psi,\psi} \leq -\underbrace{c_p D t \sep{(\langle v\rangle^{2\gamma}\lambda_v^{4s})^{\rm Wick} \psi,\psi}}_{III}.
 $$
Concerning  the third term in  \eqref{D},  let us compute $ \{p,v\cdot \xi\}$:
 \begin{align*}
 \{p,v\cdot \xi\}&= \nabla_\eta p \cdot \nabla_v (v\cdot \xi) - \nabla_v p\cdot \nabla_\eta(v\cdot \xi)= \langle v \rangle^\gamma (\nabla_\eta \lambda_v^{2s}) \cdot  \xi\\& = 2s \langle v  \rangle^\gamma  \lambda_v^{2s-2} (\eta\cdot\xi +(v\wedge \eta)\cdot(v\wedge \xi))\\
 &\leq 2 s  \langle v  \rangle^\gamma \lambda_x \lambda_v^{2s-1}
 ,\end{align*}
 where we used  the fact that $|\eta\cdot\xi +(v\wedge \eta)\cdot(v\wedge \xi)|\leq  \lambda_x \lambda_v$.
Hence,  for any $\varepsilon_{1}>0$, we obtain two non-negative terms
\begin{align*}
&Dt\sep {\{p,v\cdot \xi\}^{\rm Wick} \psi,\psi} \\
&\qquad \leq \underbrace{ 2s {\varepsilon_{1}}^{-1} D\sep{(\langle v  \rangle^\gamma\lambda_v^{2s})^{\rm Wick} \psi,\psi }}_{iii} +\underbrace{2s\varepsilon_{1}^s
D t^{1+s} \sep{ (\langle v  \rangle^\gamma\lambda_x^{s+1} \lambda_v^{s-1})^{\rm Wick} \psi,\psi}}_{iv}.
\end{align*}

Let us now consider \eqref{E}. Using the fact that $\omega \leq \langle v\rangle^\gamma \lambda_x^s \lambda_v^s $, we can bound the first term in \eqref{E}, for any $\varepsilon_{2}>0$, with two non-negative terms
 $$
E t^s \sep{\omega^{\rm Wick} \psi,\psi}  \leq \underbrace{\varepsilon_{2}^{-1}E \sep{ (\langle v\rangle^\gamma \lambda_v^{2s} )^{\rm Wick}\psi,\psi}}_{v}+ \underbrace{\varepsilon_{2} ^{1 /s} E t^{1+s} \sep{(\langle v\rangle^\gamma \lambda_x^{s+1} \lambda_v^{s-1} )^{\rm Wick}\psi,\psi}}_{vi}.
 $$
For the second term in \eqref{E},   Lemma \ref{lemma3bis} implies
 $$
 \sep{\mathcal {\it \Omega} \psi,\psi}\leq c_\omega  \sep{ (p^{3/2} q^{1/2})^\wick \psi, \psi}
 $$
and, for any $\varepsilon_{3}>0$, we have
 $$
 t^{1+s}p^{3/2} q^{1/2}\leq \varepsilon_{3}^{-1} t p^2 +\varepsilon_{3}t^{1+2s}pq.
 $$
Therefore, we can bound  the second term in \eqref{E}, using Lemma~\ref{lemma4}, for any $\varepsilon_{3}>0$, by
\begin{align*}
 &-Et ^{s+1}\sep{{\it \Omega} \psi,\psi}\\
 &\qquad \leq  c_\omega\varepsilon_{3}^{-1}E t  \sep{(p^2)^{\rm Wick} \psi,\psi }+ c_\omega\varepsilon_{3}E  t^{1+2s} \sep{ (pq)^{\rm Wick}\psi,\psi }\\
 &\qquad \leq \underbrace{2(1+K^2)c_\omega\varepsilon_{3}^{-1}E t  \sep{(\langle v \rangle^{2\gamma} \lambda_v^{4s})^{\rm Wick} \psi,\psi }}_{vii}\\
 &\qquad \quad + \underbrace{(1+K)^2c_\omega\varepsilon_{3}E  t^{1+2s} \sep{ (\langle v \rangle^{2\gamma} \lambda_v^{2s} \lambda_x^{2s})^{\rm Wick}\psi,\psi }}_{viii}
 \end{align*}
where $(vii)$ and $(viii)$ are non-negative.

\noindent Let us now observe that
$$
(\nabla_\eta \lambda_v^2)\cdot \xi=2 (\eta\cdot \xi+ (v\wedge \eta) \cdot( v\wedge \xi)),
$$
and
$$
\nabla_\eta(\eta\cdot \xi
 + (v\wedge \eta) \cdot (v\wedge \xi))\cdot \xi= \lambda_x^2- \langle v\rangle^2.
 $$
 We then compute
 \begin{align*}
&  \{\omega,v\cdot \xi\} \\
 &\quad=  \nabla_\eta \omega \cdot \nabla_v (v\cdot \xi) - \nabla_v \omega\cdot \nabla_\eta(v\cdot \xi)= \nabla_\eta \omega \cdot  \xi\\
 &\quad = - \langle v\rangle^{\gamma}\lambda_x^{s-1}\lambda_v^{s-1}\nabla_\eta(\eta\cdot \xi
 + (v\wedge \eta) \cdot (v\wedge \xi))\cdot \xi\\
 &\qquad -\langle v\rangle^{\gamma}\lambda_x^{s-1}(\eta\cdot \xi+ (v\wedge \eta) \cdot (v\wedge \xi))(\nabla_\eta\lambda_v^{s-1})\cdot \xi \\
 &\quad = - \langle v\rangle^{\gamma}\lambda_x^{s+1}\lambda_v^{s-1}
 +\langle v\rangle^{\gamma+2}\lambda_x^{s-1}\lambda_v^{s-1}
 - (s-1)\langle v\rangle^{\gamma}\lambda_x^{s-1}\lambda_v^{s-3}(\eta\cdot \xi + (v\wedge\eta)\cdot( v\wedge\xi))^2.
\end{align*}
In the last expression of $\{\omega,v\cdot \xi\}$, we first notice that since $s-1<0$ and $\min(\lambda_x,\lambda_v) \ge \langle v \rangle$, the second term is bounded as follows:
$$
\langle v\rangle^{\gamma+2}\lambda_x^{s-1}\lambda_v^{s-1} \le \langle v \rangle^{\gamma+2s}.
$$
Gathering the first and third terms, we use Cauchy-Schwarz inequality and  $s<1$ to find:
\begin{align*}
&- \langle v\rangle^{\gamma}\lambda_x^{s+1}\lambda_v^{s-1}
-(s-1)\langle v\rangle^{\gamma}\lambda_x^{s-1}\lambda_v^{s-3}(\eta\cdot \xi + (v\wedge\eta)\cdot( v\wedge\xi))^2 \\
&\quad \le - \langle v\rangle^{\gamma}\lambda_x^{s+1}\lambda_v^{s-1}
+ (1-s) \langle v\rangle^{\gamma}\lambda_x^{s-1}\lambda_v^{s-3} (\lambda_x^2-\langle v \rangle^2) (|\eta|^2 +|v \wedge \eta|^2) \\
&\quad =  - \langle v\rangle^{\gamma}\lambda_x^{s+1}\lambda_v^{s-1} + (1-s)  \langle v\rangle^{\gamma}\lambda_x^{s+1}\lambda_v^{s-3} (|\eta|^2 +|v \wedge \eta|^2) \\
&\qquad- (1-s) \langle v\rangle^{\gamma+2}\lambda_x^{s-1}\lambda_v^{s-3}(|\eta|^2 +|v \wedge \eta|^2) \\
&\quad \le  - \langle v\rangle^{\gamma}\lambda_x^{s+1}\lambda_v^{s-1} + (1-s)  \langle v\rangle^{\gamma}\lambda_x^{s+1}\lambda_v^{s-1} -
(1-s)  \langle v\rangle^{\gamma+2}\lambda_x^{s+1}\lambda_v^{s-3} \\
&\quad \le -s  \langle v\rangle^{\gamma}\lambda_x^{s+1}\lambda_v^{s-1}.
\end{align*}
Thus we have:
$$
 \{\omega,v\cdot \xi\}\le -s \langle v\rangle^{\gamma}\lambda_x^{s+1}\lambda_v^{s-1}+ \langle v\rangle^{\gamma+2s  }.
$$
Hence, the third term in  \eqref{E} can be estimated as
\begin{align*}
&E t^{s+1}\sep{\{\omega,v\cdot \xi\}^{\rm Wick} \psi,\psi} \\
&\quad \leq -\underbrace{ s Et^{s+1} \sep{(\langle v\rangle^{\gamma}\lambda_x^{s+1}\lambda_v^{s-1})^{\rm Wick} \psi,\psi})}_{IV}+\underbrace{ E t^{s+1} \sep{(\langle v\rangle^{\gamma+2s} )^{\rm Wick}\psi,\psi}}_{ix},
\end{align*}
  where $(-IV)$ is non-positive and $(ix)$ is non-negative.

  It remains to consider \eqref{F}. Observing that, for any $\varepsilon_4>0$,
$$
t^{2s}\langle v\rangle^\gamma \lambda_x^{2s}\leq \varepsilon_4 ^{-1} \langle v\rangle^\gamma \lambda_v^{2s} + \varepsilon_4^{\frac{1-s}{2s}}t^{1+s}\langle v\rangle^\gamma \lambda_v^{s-1}\lambda_x^{s+1},
$$
we have that the first term in \eqref{F} can be bounded for any $\varepsilon_4>0$, by
 \begin{align*}
 &(1+2s) t^{2s}\sep{q^{\rm Wick} \psi,\psi} \\
 &\quad \leq \underbrace{(1+2s)\varepsilon_4 ^{-1} \sep{(\langle v\rangle^\gamma \lambda_v^{2s})^{\rm Wick} \psi,\psi }}_{x} +\underbrace{(1+2s)\varepsilon_4^{\frac{1-s}{2s}}  t^{1+s}\sep{(\langle v\rangle^\gamma \lambda_v^{s-1}\lambda_x^{s+1}) ^{\rm Wick}\psi,\psi}}_{xi}\\
 &\qquad + \underbrace{K (1+2s)t^{2s}\sep{(\langle v\rangle ^{\gamma+2s})^{\rm Wick} \psi,\psi}}_{xii}
 \end{align*}
  where $(x),(xi),(xii)$ are non-negative terms.

 \noindent Moreover, using Lemma  \ref{lemma3} and Lemma \ref{lemma4}, the second term in \eqref{F} can be estimated as
 \begin{align*}
-t^{1+2s}\sep{\mathcal Q \psi,\psi}&\leq -c_q t^{1+2s}\sep{(pq)^{\rm Wick} \psi,\psi}
\le - \underbrace{c_q t^{1+2s}\sep{(\langle v \rangle^{2\gamma} \lambda_v^{2s} \lambda_x^{2s})^{\rm Wick} \psi,\psi}}_{V}
 \end{align*}
 where $(-V)$ is non-positive.
Finally, since $q$ does not depend on~$\eta$, the Poisson bracket $\{q,v\cdot \xi\}$ vanishes, hence the third term in \eqref{F} is null.

We conclude the proof as we did for Theorem~\ref{FKEthm}, checking that we can choose (in order  of reverse appearance) the constants $C$, $D$, $E$ and the small constants
$\eps_{j}$, $j=1,\dots,4$ such that for $t \in (0,1]$,
 \begin{align*}
 - I + i+iii+v+ x &\leq -\frac{1}{10}I, \\
 -II + ii+ ix+ xii 
                    &\leq -\frac{1}{10} II, \\
 -III + vii &\leq -\frac{1}{10}III,\\
 -IV + iv+vi+xi  &\leq -\frac{1}{10}IV,\\
 -V +viii &\leq -\frac{1}{10}V.
 \end{align*}
Note that $D$ and  $C$  can be taken arbitrarily larger at the end of this procedure. 
This ends the proof.
\fin

\subsection{Proof of Proposition \ref{prop:lambda1*}}

We can now prove Proposition \ref{prop:lambda1*}. Consider $\varphi$ the solution of
$$
\partial_t \varphi = v \cdot \nabla_x \varphi - A \varphi,
$$
with initial data $\varphi_0$
and $\psi=\mathcal{F}_x \varphi$ to be the  solution of
$$
\partial_t \psi { -} i v\cdot \xi \psi + A \psi =0
$$
 with initial data $\psi_0= \mathcal{F}_x \varphi_0$. From Lemma \ref{derivB}, we know that
 $$
 \hhh(t) \leq \hhh(0) = C \norm{\psi_0}^2,
 $$
 and using Lemma \ref{hh1bis}, this gives for all $t\in (0,1]$
 \begin{equation} \label{majorpq}
 \sep{p ^{\rm Wick} \psi,\psi} \leq \frac{2C}{D} \frac{1}{t} \norm{\psi_0}^2 \qquad \textrm{ and }  \qquad \sep{q^{\rm Wick} \psi,\psi} \leq  \frac{2C}{t^{1+2s}} \norm{\psi_0}^2,
 \end{equation}
 where we used the fact that both left members are non-negative according to Proposition~\ref{mainpseudo}. Working in the class $S_K(p)$ again, gives through Proposition~\ref{mainpseudo} and Lemma~\ref{AHLinverse} (see there the definition of $H_R$)
 \begin{equation*}
 \begin{split}
 \norm{ \seq{v}^{\gamma/2} \seq{D_v}^s \psi}^2
& =  \norm{ \seq{v}^{\gamma/2} \seq{D_v}^s ({ (p^{1/2})}^w)^{-1} { (p^{1/2})}^w \psi}^2 \\
& = \Vert \underbrace{\seq{v}^{\gamma/2} \seq{D_v}^s ({ (p^{1/2})}^{-1})^w}_{\textrm{bounded operator}} H_R \, { (p^{1/2})}^w \psi \Vert^2 \\
& \lesssim \norm{  { (p^{1/2})}^w \psi}^2, \\
\end{split}
\end{equation*}
where we used that the  operator $ \seq{v}^{\gamma/2} \seq{D_v}^s $ has its Weyl symbol in $S_K(p^{1/2})$ (this Weyl symbol is $\seq{v}^{\gamma/2} \sharp \seq{\eta}^s$), and that $ {(p^{1/2})}^{-1} \in S_K(p^{-1/2})$ , so that $\seq{v}^{\gamma/2} \seq{D_v}^s ({ (p^{1/2})}^{-1})^w$ is a bounded operator.
Using then \eqref{mainnorm} and  \eqref{majorpq}, we get
\begin{equation*}
 \begin{split}
& \norm{ \seq{v}^{\gamma/2} \seq{D_v}^s \psi}^2  \lesssim \norm{  { (p^{1/2})}^w \psi}^2
 \simeq \sep {p^\wick \psi,\psi}
 \lesssim \frac{1}{t} \norm{\psi_0}^2.
\end{split}
\end{equation*}
Similarly,
$$
\norm{ \seq{v}^{\gamma/2+s} \psi}^2  \lesssim \frac{1}{t} \norm{\psi_0}^2,
$$
and working in $S_K(q)$ gives, in the same way,
$$
\norm{ \seq{v}^{\gamma/2}\seq{\xi}^s \psi}^2  \lesssim \frac{1}{t^{1+2s}} \norm{\psi_0}^2.
$$
Taking the inverse Fourier transform in the $x$ variable finally yields
\begin{multline*}
\norm{ \seq{v}^{\gamma/2} \seq{D_v}^s \varphi}^2 \lesssim \frac{1}{t} \norm{\varphi_0}^2,  \qquad \norm{ \seq{v}^{\gamma/2+s} \varphi}^2  \lesssim \frac{1}{t} \norm{\varphi_0}^2 \\
\textrm{ and } \qquad \norm{ \seq{v}^{\gamma/2}\seq{D_x}^s \varphi}^2  \lesssim \frac{1}{t^{1+2s}} \norm{\varphi_0}^2.
\end{multline*}
This is exactly the statement of Proposition \ref{prop:lambda1*}, the proof is thus complete. \fin

\section{Adaptation of the proof for the primal result and generalization} \label{sec:general}

{\subsection{Adaptation of the proof for the primal result} \label{adaptprimal}
If we want to prove the ``primal'' regularization property in Theorem~\ref{thm:LIBthm}, as in Subsection~\ref{subsec:splitting}, we split $\Lambda$ into two parts:
\begin{align} \label{decomptilde}
\Lambda f &=\left(- K \langle v \rangle^{\gamma+2s} - v \cdot \nabla_x f + \int_{\R^3 \times \mathbb{S}^2} B_\delta(v-v_*,\sigma) \mu'_* (f'-f) \right)\, \d\sigma \, \dv_* \\
&\qquad + \Bigg(K \langle v \rangle^{\gamma+2s} + \int_{\R^3 \times \mathbb{S}^2} B_\delta(v-v_*,\sigma) (\mu'_*-\mu_*) f\, \d\sigma \, \dv_*\\
&\hskip 4.5cm + \int_{\R^3\times \mathbb{S}^2} B_\delta^c(v-v_*,\sigma) (\mu'_* f' - \mu_* f) \, \d\sigma \, \dv_* + Q(f,\mu) \Bigg) \\
&=: \widetilde{\Lambda}_1f + \widetilde{\Lambda}_2f.
\end{align}
Then, the study of $\widetilde{\Lambda}_1^m$ is totally similar to the one of $\Lambda_1^{*,m}$ (the only differences being in the fact that the roles of $m$ and $m^{-1}$ are inverted and the sign in front of the transport operator is opposite). We thus just have to adapt the signs in the Lyapunov functional: the sign of $\omega$ has to be changed in Subsection~\ref{subsec:refweights}. The other part $\widetilde{\Lambda}_2$ is controlled as well as $\Lambda_2$. The proof is thus done in the same way and we do not enter into details.}

\subsection{Generalization to higher order estimates}
Theorem \ref{thm:LIBthm} deals with regularization in close to $L^2$ spaces: for example, it says that that the semigroup associated to $\Lambda=L-v \cdot \nabla_x$ with $L$ given in \eqref{eq:BO}
goes from $L^2$ to $H^{s}$ type spaces, with suitable weights and explicit norms. One can wonder if an higher order quantitative regularization is also available. This is the aim of the following Theorem, for which we give a condensed statement in the primal case and in  homogeneous $H^{\ell s}$ spaces  (see notation \eqref{homonorm} and below).

\begin{thm} \label{thm:LIBthmhigh}
 Let $\ell\in \N^*$ and  {$k' \ge 0$, $k>\max(\gamma/2+3+2\ell s, k' + \gamma+5/2)$}. Let  $f$ be a solution of \eqref{FKE} with $L$ given by~\eqref{eq:BO} with initial data  $f_0 \in H^{\ell s}_{x,v}(\langle v \rangle^k)$. Then, there exists a constant $C_\ell>0$ independent of~$f_0$ such that we have the following regularization estimates: for any  $t \in (0,1]$ we have
 $$
  \|f(t)\|_{H^{\ell s}_{x,v}(\langle v \rangle^k)} \le \frac{C_{\ell}}{t^{1/2+s}} \|f_0\|_{H^{(\ell-1)s}_{x,v}(\langle v \rangle^{k'})}.
$$
%
\end{thm}

In this Section we shall not give the complete proof of this result, since this is very close to the  one of Theorem \ref{thm:LIBthm}, but only elements of it. The remaining of this Section is devoted to these elements.

\bigskip
\noindent {\bf Elements of proof of Theorem \ref{thm:LIBthmhigh}.} \ \
In all the following, we consider $\ell \in \N^*$ given by the theorem as well as $k$ and $k'$ given there. Recall that the statement gives a (primal) regularization result on the solution $f(t)$ of $\D_t f = \Lambda f $ where $\Lambda = -v \cdot \nabla_x + L$ where $L$ is the linearized Boltzmann collision kernel given  in \eqref{eq:BO}.

As a first step we  split operator $\Lambda$ into two parts following  \eqref{decomptilde}
$$
\Lambda = \tilde{\Lambda}_1 + \tilde{\Lambda}_2.
$$
Adapting the proof of Lemma \ref{lem:lambda_2}, we have for suitable functions $h$
\begin{equation} \label{eq:lambda_2bis}
\|\tilde{\Lambda}_2  h \|_{H^{(\ell-1) s}_{x,v}( \seq{v}^k)} \lesssim \|h\|_{H^{(\ell-1) s}_{x,v} (\langle v \rangle^{k'} )}
\end{equation}
where $k$ and $k'$ are given in the statement of Theorem \ref{thm:LIBthmhigh}.
We shall in a moment  prove that

\begin{prop} \label{resultHs}
We have for all $t \in (0,1]$,
$$
\norm{S_{\tilde{\Lambda}_1}(t) h}_{H^{\ell s}_{x,v}(\seq{v}^k)} \lesssim
\frac{1}{t^{1/2+s}} \|h\|_{H^{(\ell-1)s}_{x,v}(\langle v \rangle^{k'})}.
$$
\end{prop}

Taking this result into account and together with \eqref{eq:lambda_2bis}  we can write
$$
 S_{\Lambda}(t) = S_{\tilde{\Lambda}_1}(t) + \int_0^t S_\Lambda(t-s)  (\tilde{\Lambda}_2 S_{\tilde{\Lambda}_1})  (s) \, ds
 $$
  for $t\in [0,1)$. Arguing as in the proof of Lemma \ref{lem:lambda_1/lambda} we easily get the Theorem  (this strongly uses $s<1/2$). We omit the details.
  \fin

\bigskip
\noindent {\bf Elements of proof of Proposition \ref{resultHs}.} \ \
 We notice that it is sufficient to prove the following two estimates :
\begin{equation} \label{FGHl}
\norm{S_{\tilde{\Lambda}_1^{m^{-1}}(t) } \phi}_{G_\ell} \lesssim  \frac{1}{t^{1/2+s}} \norm{\phi}_{F_{\ell-1}}, \qquad \norm{S_{\tilde{\Lambda}_1^{m^{-1}}(t) } \phi}_{H_\ell} \lesssim  \frac{1}{t^{1/2}} \norm{\phi}_{F_{\ell-1}}
\end{equation}
where similarly to the beginning of Subsection \ref{subsec:reg}, we define
(here in the primal case)
$$
\left\{
\begin{aligned}
&F=L^2_{x,v} \\
&G_\ell =H^{\ell s,0}_{x,v}({ \langle v \rangle^{\ell \gamma/2}}) \\
&H_\ell=H^{0,\ell s}_{x,v}( \langle v \rangle^{\ell \gamma/2}) \cap  L^2_{x,v}( \langle v \rangle^{\ell (\gamma+2s)/2})
\end{aligned}
\right.
$$
and $\tilde{\Lambda}_1^{m^{-1}} = m^{-1} \tilde{\Lambda}_1 m$.
The proof is very close to the one given in the dual case in the Section 3. As mentioned in the previous subsection, we essentially have to replace $m$ there by $m^{-1}$ here, change the sign in front of the drift $v.\nabla_x$, we also have to work in $G_\ell$ or $H_\ell$ instead of $G (= G_1)$ and $H (=H_1)$ introduced in Subsection \ref{subsec:reg} for getting Proposition \ref{prop:lambda1*}.

In fact by interpolation, estimates \eqref{FGHl} are direct consequences of the following estimates:
\begin{equation} \label{FGHl2}
\norm{S_{\tilde{\Lambda}_1^{m^{-1}}(t) } \phi}_{G_\ell} \lesssim  \frac{1}{t^{\ell(1/2+s)}} \norm{\phi}_{F}, \qquad \norm{S_{\tilde{\Lambda}_1^{m^{-1}}(t) } \phi}_{H_\ell} \lesssim  \frac{1}{t^{{\ell}/2}} \norm{\phi}_{F}.
\end{equation}

We shall in fact give an idea on how to prove the preceding result using the same tools as in Section \ref{sec:quantization}. Let us recall that a fundamental large parameter $K$ is involved there and enters here in the definition of $\tilde{\Lambda}_1^{m^{-1}}$. Following the strategy of Section
\ref{sec:quantization}, we get that
$$
\tilde{\Lambda}_1^{m^{-1}} = -b^w - v\cdot \nabla_x
$$
where $b$ has exactly the same properties than $a$ in Section~\ref{sec:quantization}. In particular as in Lemma~\ref{aposell}, $\Re b \geq 0$ and $\Re(b)$ is elliptic positive in the class $S_K(p)$ as there.
We pose $B = b^w$  and recall the definitions of the symbols in Subsection \ref{subsec:refweights}: for given $s \in (0,1/2)$ and $\gamma \in (0,1)$,
$$
p(v,\eta) = \seq{v}^\gamma \lambda_v^{2s} + K \seq{v}^{\gamma + 2s},
$$

$$
q(v,\eta) = \seq{v}^\gamma \lambda_x^{2s} + K \seq{v}^{\gamma + 2s},
$$
and
$$
\omega (v,\eta) = { - }\langle v\rangle ^\gamma \lambda_x^{s-1} \lambda_v^{s-1}(\eta\cdot \xi+ (v\wedge \eta) \cdot (v\wedge \xi)).
$$
Since we are in the primal and not dual case (he sign in front of the transport term is opposite), we have to take the opposite of $\omega$ that we call $\tilde{\omega}:=-\omega$.

The main point of the analysis is then to introduce, such as in Subsection~\ref{subsec:proplambda1*}, a suitable functional which is here:
\begin{multline} \label{hhhl}
\hhh_\ell(t) := C\norm{\psi}^2 +
\sum_{0 \leq \alpha+ \beta \leq  \ell-1}
D_{\alpha,\beta} t ^{1+\alpha + \beta (1+2s)}\sep{\sep{p^{1+\alpha}q^{\beta}  }^{\rm Wick} \psi,\psi} \\
+ E_{\alpha,\beta} t^{1/2+ \alpha + (1/2+\beta)(1+2s)}  \sep{ \sep{p^{\alpha}q^{\beta} \tilde{\omega}} ^{\rm Wick} \psi,  \psi} \\
+ F_{\alpha,\beta}t^{\alpha + (1+\beta)(1+2s)} \sep{\sep{p^{\alpha}q^{1+\beta}}^{\rm Wick} \psi,\psi}
\end{multline}
 for well chosen constants $C$, $D_{\alpha,\beta}$, $ E_{\alpha,\beta}$ and $F_{\alpha,\beta}$.
 We note that for $\ell=1$, we get $\hhh_1 = \hhh$ defined in \eqref{hhh1}. The computations exactly follow the ones done in Subsection~\ref{subsec:proplambda1*} using estimates
similar to the ones given in Subsection~\ref{subsec:technical}, with the same roles of each term as there in the preceding decomposition. Note that we were note able to restrict the analysis to $\alpha+\beta = \ell-1$ due to too high order terms after time derivation, this explains that the full range of $\alpha$ and $\beta$ is needed to close the estimates and conclude that
$$
\frac{d}{ dt} \hhh_l(t) \leq 0.
$$
We omit the details of the computation as well as the last parts of the proof of \eqref{FGHl2} which leads to Proposition \ref{resultHs}, since it follows the end of Section \ref{sec:quantization} . \fin

\Bk

\appendix \section{}
\subsection{Carleman representation}
We  state here a classical tool in the analysis of Boltzmann operator: the Carleman representation. We refer to~\cite{AHL*} for more details on the version that we state here.
\begin{lem}[Carleman representation] \label{lem:Carleman}
	Let $F$ be a measurable function defined on~$(\R^3)^4$. For any vector $h \in \R^3$, we denote by $E_{0,h}$ the (hyper)vector plane orthogonal to $h$. Then, when all sides are well defined, we have the following equality :
		\begin{align*}
		&\int_{\R^3 \times \S^2} b(\cos \theta) |v-v_*|^\gamma F(v,v_*,v',v'_*) \, \dv_* \, \d\sigma \\
		&\quad = \int_{\R^3_k} \d h \int_{E_{0,h}} \d\alpha \, \tilde{b}(\alpha,h)\,\mathds{1}_{|\alpha| \ge |h|} \frac{|\alpha+h|^{\gamma+1+2s} }{|h|^{3+2s}} \,F(v,v+\alpha-h,v-h,v+\alpha)
		\end{align*}
	where $\tilde{b}(\alpha,h)$ is bounded from above and below by positive constants and $\tilde{b}(\alpha,h)=\tilde{b}(\pm\alpha,\pm h)$.
\end{lem}

\subsection{Pseudodifferential calculus} \label{pseud}

We first recall the definitions of the quantizations we shall use in the following.
Let us consider a temperate symbol $\sigma \in \sss$, we define its standard quantization $\sigma(v,D_v)$ for  $f\in L^2(\R^d)$ by
$$
\sigma(v,D_v) f (v) := \frac{1}{(2\pi)^d} \int e^{iv \cdot\eta} \sigma(v,\eta) \hat{f}(\eta) \, \d \eta.
$$
The Weyl quantization is defined by
$$
\sigma^w f (v) := \frac{1}{(2\pi)^d} \iint e^{i(v-w) \cdot\eta} \sigma\left(\frac{v+w}{2},\eta\right) f(w) \, \d \eta \, \d w.
$$
We recall that  for two symbols $\sigma$ and $\tau$ we have
\begin{equation} \label{defsharp}
{\sigma^w \tau^w = (\sigma \sharp \tau)^w, \quad \sigma \sharp \tau = \sigma\tau + \int_0^1 (\D_\eta \sigma \sharp_\theta \D_v \tau - \D_v \sigma \sharp_\theta \D_\eta \tau ) \, \d \theta}
\end{equation}
where for $V=(v,\eta)$ we have $\sharp = \sharp_1$ and for $\theta \in (0,1]$,
$$
\sigma \sharp_\theta \tau (V):= \frac{1}{2i} \iint e^{-2i [V-V_1,V-V_2]/\theta} \sigma(V_1) \tau(V_2)\, \d V_1 \, \d V_2 /(\pi \theta)^d
$$
with $[V_1,V_2] = v_2\cdot\eta_1 - v_1\cdot\eta_2$ the canonical symplectic form on $\R^{2d}$.
We shall also use the Wick quantization, which has very nice properties concerning positivity of operators (see~\cite{Lerner-CME,LernerBook1,LernerBook2} for more details on the subject). For this, we first introduce the Gaussian in  phase variables
\begin{equation} \label{wickN}
N(v,\eta) := (2\pi)^{-d} \exp^{-(|v|^2+|\eta|^2)/2}.
\end{equation}
The Wick quantization is then defined by
\begin{equation} \label{wickN2}
\sigma^\wick f(v) := (\sigma \star N )^w f(v),
\end{equation}
where $\star$ denotes the usual convolution in $(v,\eta)$ variables.
Recall that one of the main property of Wick quantization is its positivity:
\begin{equation} \label{wick2}
\forall \, (v,\eta) \in \R^6, \; \sigma(v,\eta)\geq 0 \Rightarrow \sigma^{\rm Wick}\geq 0,
\end{equation}
and that the following relation holds (see e.g. \cite[Proposition 3.4]{Lerner-CME}):
\begin{equation}\label{wick}
[ g^{\rm Wick},iv\cdot \xi]=\{  g, v\cdot \xi\}^{\rm Wick}.
\end{equation}
The previous definitions extend to symbols in $\sss'$ by duality.

\subsection{The weak semiclassical class $S_K(g)$}

 Let $\Gamma := |\dv|^2 + |\d \eta|^2$ be the flat metric on~$\R^6_{v,\eta}$. The first point is to verify that the introduced symbols and weights are indeed   in a suitable symbolic calculus with large parameter $K$ uniformly in the parameter $\xi$.
For this, we first recall  that a weight $1 \leq g$ is said to be temperate with respect to $\Gamma$ if there exist $N \geq 1$ and $C_{N}$ such that  for all $(v, \eta)$, $(v',\eta') \in \R^6$
$$
g(v',\eta') \leq C_N \, g(v,\eta) ( 1 + |v'-v| + |\eta'-\eta|)^N
$$
We now introduce adapted classes of symbols.

\begin{defn}
Let  $g$ be a temperate weight.  We denote by $S(g)$ the symbol class of all smooth functions $\sigma(v,\eta)$ (possibly depending on parameters $K$ and $\xi$) such that
$$
\abs{\D_v^\alpha \D_\eta^\beta \sigma(v,\eta) } \leq C_{\alpha,\beta} g(v,\eta)
$$
where for any multiindex $\alpha$ and $\beta$, $C_{\alpha,\beta}$ is uniform in $K$ and $\xi$. We denote also $S_K(g)$ the symbol class of all smooth functions $\sigma(v,\eta)$ (possibly depending on $K$ and $\xi$ again) such that
$$
\abs{\sigma} \leq C_{0,0} g \quad \text{and} \quad \forall \, |\beta| \geq 1, \quad \abs{\D_v^\alpha \D_\eta^\beta \sigma } \leq
 C_{\alpha,\beta} K^{-1/2} g
$$
uniformly in $K$ and $\xi$. Note that $S_K(g) \subset S(g)$ and that these definitions are with respect to the flat metric.

Eventually, we shall say that a symbol $\sigma$ is elliptic positive in $S(g)$ or $S_K(g)$ if in addition $\sigma \geq 1$ and there exists a constant $C$ uniform in parameters such that $C^{-1} g \leq \sigma \leq C g$.
\end{defn}
Before focusing on the class $S_K(g)$, we first recall one of the main results concerning the class without parameter (and without weight) $S(1)$:

\begin{lem}[Calderon Vaillancourt Theorem]
Let $\sigma \in S(1)$, then $\sigma^w$ is a bounded operator with norm depending only on a finite number of semi-norms of $\sigma$ in $S(1)$.
\end{lem}
The classes $S_K$ and $S$ have standard internal properties:

\begin{lem} \label{stab} For $K$ sufficiently large, we have the following:
\begin{enumerate}[label=\alph*)]
\item Let $g$ be a temperate weight and consider $\sigma$ an elliptic positive symbol in $S_K(g)$ then for all $\nu \in \R$, $\sigma^\nu \in S_K(g^\nu)$;
\item Let $g$, $h$ be temperate weights and consider $\sigma$ in $S_K(g)$, $\tau$ in $S_K(h)$, then $\sigma\tau$ is in $ S_K(g h)$.
\end{enumerate}
\end{lem}

\preuve For point \it a)\rm, just notice that if $\sigma$ is an elliptic positive symbol in $S_K(g)$, then $\sigma \simeq g$ so that $\sigma^\nu \simeq g^\nu$. We also have directly for $\beta$ a multiindex of length 1
$$
\abs{ \D_\eta^\beta \sigma^\nu } = |\nu| \sigma^{\nu-1} \abs{\D_\eta^\beta \sigma} \leq C g^{\nu-1} K^{-1/2} g = C K^{-1/2} g^\nu
$$
using $\sigma \simeq g$. Estimates on  higher order derivatives are straightforward.

For point \it b)\rm, the computation is also straightforward using the Leibniz rule.
\fin

Now we can quantize the previously introduced symbols. The main \it semiclassical \rm idea behind the introduction of the class $S_K$ for $K$ large is that invertibility and powers of operators associated to symbols are direct consequences of similar properties of symbols, essentially independently of the quantization.

We first check that the class $S_K$ is essentially stable by change of quantization.

\begin{lem} \label{allsk}
Let $g$ be a temperate weight and consider $\tilde{\sigma}$  a positive elliptic symbol in $S_K(g)$. We denote ${\sigma}$ the Weyl symbol of the operator  $\tilde{\sigma}(v,D_v)$ so that
$\sigma^w = \tilde{\sigma}(v,D_v)$  and recall that the Weyl symbol of $\sigma^\wick$ is $\sigma \star N$.  Then
 $\sigma$ and  $ \sigma\star N$ are both in $S_K(g)$. If in addition $\tilde{\sigma}$ is elliptic positive, then $\Re \sigma$ and $\Re \sigma\star N$ are elliptic positive.
\end{lem}

\preuve
We first prove the result for $\sigma$ supposing that $\tilde{\sigma}$ is elliptic positive. From for e.g.~\cite{LernerBook2} and an adaptation of Lemma 4.4 in~\cite{AHL*}, we know that
\begin{equation} \label{diff}
\sigma - \tilde{\sigma} \in K^{-1/2} S(g).
\end{equation}
Since $K^{-1/2} S(g) \subset S_K(g)$, this gives that ${\sigma} \in S_K(g)$.
If in addition $\tilde{\sigma}$ is elliptic positive, then let us prove that   $\Re \sigma$ also is.  There exist constants $C$, $C'$ uniform in $K$ large such that
$$
C^{-1} g - C'K^{-1/2} g \leq \Re {\sigma} \leq  C g + C'K^{-1/2} g
$$
if $C^{-1} g \leq \sigma  \leq  C g$. Taking $K$ sufficiently large then gives the result.

We now deal with $\sigma\star N$, supposing that $\sigma$  is in $S_K(g)$. For $V=(v,\eta)$ we have
$$
\sigma\star N (V)
 = \iint \sigma(V-W) N(W) \d W
 $$
 and using the temperance property of $g$, we get uniformly in all other possible parameters (including $K$)
 $$
 \abs{\sigma\star N(V)} \leq \iint C g(V) (1+|W|)^N N(W) \, \d W \leq C' g(V).
 $$
 For the derivatives, we get similarly for multiindex $\alpha$ and $\beta$ with $|\beta| \geq 1$
 \begin{equation}
 \begin{split}
 \abs{\D^\alpha_v \D^\beta_\eta \sigma\star N (V)}
 &  \leq  \iint \abs{\D^\alpha_v \D^\beta_\eta \sigma(V-W)} N(W) \, \d W\\
 & \leq C K^{-1/2} \iint g(V-W) N(W) \, \d W \\
 & \leq C' K^{-1/2} \iint g(V) (1+|W|)^N N(W) \, \d W \\
 & \leq C'' K^{-1/2} g(V).
 \end{split}
 \end{equation}
Suppose now that in addition
 $\tilde{\sigma}$  is elliptic positive, then $\Re \sigma$ is {elliptic positive} and $C^{-1} g(V) \leq \Re \sigma(V) \leq C g(V)$ for a constant $C>0$. Since $\Re \sigma\star N$ is positive, this implies with the temperance of $g$ that
\begin{multline}
c' g(V) { \le} \iint C^{-1} C_N^{-1} g(V) (1+|W|)^{-N} N(W) \d W \\ \leq \Re \sigma\star N(V)  \leq \iint C C_N g(V) (1+|W|)^N N(W) \d W = C' g(V)
\end{multline}
for some positive constants $c'$ and $C'$, so that $\Re \sigma\star N$ is indeed elliptic positive.
 \fin

\remark \label{allskrem} Note that using exactly the same argument as in the proof before, we also get that if $\tau$ is a given  elliptic positive symbol in $S_K(g)$, with $g$ a temperate weight,  then $\tau \star N$  is also an elliptic positive symbol in $S_K(g)$.

\bigskip
The next technical lemma is also proven in~\cite{AHL*}:

\begin{lem}[Lemma 4.2 in \cite{AHL*}] \label{AHLinverse} Let $g$ be a temperate weight and $\sigma \in S_K(g)$. Then for $K$ sufficiently large (depending on a finite number of semi-norms of $\sigma$), the operator $\sigma^w$ is invertible and there exists $H_L$ and $H_R$ bounded invertible operators that are close to identity as well as their inverse such that
$$
(\sigma^w)^{-1} = H_L (\sigma^{-1})^w = (\sigma^{-1})^w H_R.
$$
The norms of operators $H_L$ and
$H_R$ and their inverse can be bounded uniformly in parameters (including $K$).
\end{lem}
Note that by ``close to  identity uniformly in parameters'', we mean that
$$
\norm{H_L f} \simeq \norm{H_R f} \simeq \norm{f}.
$$
with constants uniform in parameters (including $K$ sufficiently large).

\preuve The proof follows exactly the lines of the one given in \cite[Lemma 4.2. i)]{AHL*}.
\fin

%



We now give the main Proposition that will be used in the proof of the technical Lemmas in Subsection \ref{subsec:technical}.

\begin{prop} \label{mainpseudo} Let $g$ be a temperate  weight and consider $\sigma$ an elliptic positive symbol in $S_K(g)$. Then for $K$ sufficiently large, we have the following
\begin{equation} \label{maingamma}
 \norm{(\sigma^w)^{1/2} f} \simeq  \norm{(\sigma^{1/2})^w f} \qquad \textrm{ and } \qquad
  \norm{(\sigma^w)^{-1} f} \simeq  \norm{(\sigma^{-1})^w f}.
\end{equation}
In addition, suppose that $\tau$ is another elliptic positive symbol in $S_K(g)$ then
\begin{equation} \label{mainequiv}
\norm{\sigma^w f} \simeq \norm{\tau^w f}.
\end{equation}
In particular, we have
\begin{equation} \label{mainnorm}
\norm{\sigma^w f}^2 
\simeq \norm{\sigma^\wick f}^2 \simeq \sep{ (\sigma^2)^\wick f,f}
\end{equation}
and
\begin{equation} \label{mainscalar}
\sep{\sigma^w f,f} 
\simeq \sep{\sigma^\wick f,f}
\end{equation}
uniformly in parameters (in particular $K$).
\end{prop}

\preuve
We first prove \eqref{maingamma}. For the second almost equality, we just have to notice that from Lemma \ref{AHLinverse}, we have
$$
 \norm{(\sigma^w)^{-1} f} =  \norm{H_L(\sigma^{-1})^w f} \simeq \norm{(\sigma^{-1})^w f}
 $$
 since $H_L$ is close to identity (uniformly in parameters). For the first part of~\eqref{maingamma}, we write that
\begin{equation} \label{eqnormprod}
\begin{split}
\norm{\sigma^w f}^2 & =  ( (\sigma \sharp \sigma)^w f,f) =
 ( (\sigma^2)^w f,f) + ( r^w f,f)
\end{split}
\end{equation}
where $r = \sigma\sharp \sigma - \sigma^2 \in K^{-1/2} S(g^2)$ by standard symbolic calculus. More precisely, we can write from \eqref{defsharp}
$$
r = \int_0^1 (\D_v \sigma \sharp_\theta \D_\eta\sigma  -\D_\eta\sigma \sharp_\theta \D_v \sigma) \, \d\theta
$$
and using that $\D_v \sigma \in S(g)$ and $\D_\eta \sigma \in K^{-1/2} S(g)$ gives the result by stability of the flat symbol class $S(g)$. We therefore get that
\begin{equation*}
\begin{split}
 |( r^w f,f)| & = \big| \big(  (\sigma^w)^{-1}r^w (\sigma^w)^{-1} \sigma^w f, \sigma^w f ) \big| \\
 & = \big| \big( H_L (\sigma^{-1})^w r^w (\sigma^{-1})^w H_R \sigma^w f, \sigma^w f ) \big|.
 \end{split}
 \end{equation*}
Now $\sigma^{-1} \sharp r \sharp \sigma^{-1} \in K^{-1/2} S(1)$ since $\sigma^{-1} \in S(g)$, so that $(\sigma^{-1})^w r^w (\sigma^{-1})^w$ is a bounded operator with norm controlled by a constant times $K^{-1/2}$. Since $H_L$ and $H_R$ are bounded operators independently of $K$,  there exists a constant such that
 $$
 |( r^w f,f)| \leq C K^{-1/2}\norm{\sigma^w f}^2.
 $$
This estimate and \eqref{eqnormprod}, gives that for $K$ sufficiently large,
 \begin{equation} \label{eqnormprodbis}
\begin{split}
\frac{1}{2} \norm{\sigma^w f}^2 \leq
 ( (\sigma^2)^w f,f) \leq 2\norm{\sigma^w f}^2.
\end{split}
\end{equation}
Taking $\sigma^{1/2} \in S_K(g^{1/2})$ (by Lemma \ref{stab})  instead of $\sigma$, we obtain
 \begin{equation*}
\begin{split}
 \norm{(\sigma^{1/2})^w f}^2 \simeq
 ( \sigma^w f,f) & = \norm{(\sigma^w)^{1/2}}^2
\end{split}
\end{equation*}
and the proof of \eqref{maingamma} is complete.

Concerning \eqref{mainequiv}, we just have to prove one inequality since the result is symmetric in $\tau$ and $\sigma$. For $K$ sufficiently large, we have
$$
\norm{\tau^w f} = \norm{\tau^w (\sigma^w)^{-1} \sigma^w f} = \norm{\tau^w (\sigma^{-1})^w H_R \sigma^w f} = \norm{( \tau \sharp (\sigma^{-1}))^w H_R \sigma^w f} \leq C \norm{ \sigma^w f}
$$
 since $\tau \sharp (\sigma^{-1}) \in S(1)$, so that $( \tau \sharp (\sigma^{-1}))^w$ is bounded (with bound independent of $K$). By symmetry, this proves \eqref{mainequiv}.

We then prove \eqref{mainnorm}. We first recall that $\sigma^\wick = (\sigma\star N)^w$ and that $\sigma\star N$ is elliptic positive in $S_K(g)$ by Lemma \ref{allsk}.
 From \eqref{mainequiv}, this directly yields
 $$
 \norm{ \sigma^w f} \simeq \norm{ (\sigma \star N)^w f} = \norm{ \sigma^\wick f}.
 $$
 By direct computation $(\sigma^2\star N)^{1/2}$ is also in $S_K(g)$ by
 point b) of Lemma \ref{stab} with  $\nu=2$ and $\nu=1/2$, respectively, and Lemma \ref{allsk}. Using again
 \eqref{mainequiv} and  \eqref{maingamma}, yields
 that
 $$
  \norm{ \sigma^w f} \simeq  \norm{ ((\sigma^2\star N)^{1/2})^w f} \simeq \norm{ ((\sigma^2\star N)^w)^{1/2} f} = ((\sigma^2\star N)^w f,f) = ((\sigma^2)^\wick f,f).
  $$

%
%

The proof of the last point \eqref{mainscalar} follows exactly the same lines and we skip it.
\fin


\bigskip
\bibliographystyle{acm}

\bibliography{biblioBoltzmann}

\end{document}